\documentclass[letterpaper, 10 pt, conference]{ieeeconf}
\usepackage{cite}
\usepackage{amsmath,amssymb,amsfonts}
\usepackage{algorithmic}
\usepackage{graphicx}
\usepackage{textcomp}
\usepackage{cite}
\usepackage{amsmath,amssymb,amsfonts}
\usepackage{algorithmic}
\usepackage{graphicx}
\usepackage{textcomp}

\usepackage{epsfig}
\usepackage{latexsym}
\usepackage{pdfsync}

\newcommand{\bc}{\begin{center}}
\newcommand{\ec}{\end{center}}

\newcommand{\eq}{\begin{equation}\begin{array}{rllllllllllllllllllllllllllllllll}}
\newcommand{\ee}{\end{array}\end{equation}}
\newcommand{\bmt}{\left[ \begin{array}{ccccccccc}}
\newcommand{\emt}{\end{array}\right]}
\newcommand{\bea}{\begin{eqnarray}}
\newcommand{\eea}{\end{eqnarray}}
\newcommand{\bean}{\begin{eqnarray*}}
\newcommand{\eean}{\end{eqnarray*}}

\def\tr {^{\rm T}}

\begin{document}
\title{Boundary  Stabilization of a \\Bending and  Twisting Wing\\
by Linear Quadratic Gaussian Theory
 }
\author{Arthur J. Krener,  \IEEEmembership{Life Fellow}
\thanks{ This work was funded by AFOSR under grant FA9550-23-1-0318.}
\thanks{A. J. Krener is with the Department of Mathematics, University of California, Davis, CA 95616 , USA, ajkrener@nps.edu.}}

\maketitle

\begin{abstract}
We first consider the  stabilization of the bending and twisting of a rectangular cantilever beam of  moderate to high aspect ratio using full state feedback boundary control.   Our approach is an infinite dimesnional extension of Linear Quadratic Regulation (LQR).  The we develop an infinite dimensional Kalman filter that processes two point measurements and returns an estimate of the full state.    The Linear Quadartuc Gaussian approach is to use this estimate in the place of the full state in the LQR feedback
 Then we add aerodynamic forces to obtain a model of a wing.  The aerodyamic model is based on a two dimensional state space approximate realiztion of Wagner's indicial function by R.~T.~Jones.

\end{abstract}

\section{Introduction}
 High Altitude Long Endurance  (HALE)  aircraft  have highly flexible wings.
Helicopter blades are also highly  flexible.  One would like to use active control to dampen their oscillations.   One can model a wing 
by the beam equation or the plate equation. In \cite{HMA23} it is stated that "even for relatively low aspect ratios $(4<AR<5)$ of thin rectangular cantilever plates, the model based on the
beam theory is closer to (our) experimental model (as) compared to the plate theory." So we choose a model based on the beam equation.  
   
   We assume that the beam can bend and  twist and that there are two point actuators, one can deliver a bending moment to the root of the beam and the other can deliver a torque also at  the root of the beam.  We also assume that there are two point sensors, one measures the vertical velocity of the tip of the beam and other 
measures the angular velocity of  the tip of the beam.   We wish to design a compensator that processes the measurements and delivers a bending moment and torque that dampens the bending and torsion oscillations.  We  design this compensator using  Linear Quadratic Gaussian (LQG) theory.  After accomplishing this we convert the model of the beam into a model of a wing by adding the aerodynamic forces that the motions of the beam generate.   We model the aerodynamic forces by a state space approximate realization of the classical indicial function of Wagner.  Wagner's indicial function is for an airfoil, a cross section of a wing, but we extend it to the whole wing.

Linear Quadratic Gaussian (LQG) theory is a well-known  way  of designing a dynamic compensator for a controlled and observed finite dimensional dynamical system.
We extend LQG to designing a dynamic compensator for a controlled and observed infinite dimensional dynamical system under point actuation and point sensing.  The LQG approach breaks the problem into two parts.  The first part uses Linear Quadratic Regulator theory to find a stabilizing linear feedback.   This feedback assumes that exact measuremnts of the full state are available.   In most situations this assumption is unrealistic.  The second part is to design a Kalman filter that processes the noisey measurements and returns an estimate of the full state.  This estimate is used in place of actual state in the linear feedback.  The result is a dynamic compenstor.  The closed loop eigevalues of the combined system are the union of the eigenvalues of system under full state feedbck and the eigenvalues of the error dynamics of the Kalman filter.  If all these eigevalues lie in the open left plane then the combined plant and Kalman filter is asymptotically stabe.

We start with the  model presented in Section 3.6 of the classic treatise of Bisplinghoff, Ashley and Halfman,  \cite{BAH96}.   This is a  linear model which ignores the nonlinear interactions between the bending and the  torsion of the  rectangular beam.  We add point actuation and point  sensing to the boundary of the model.  Point actuation and point  sensing are idealizations of actuation and sensing over small domains on the boundary and the idealization simplifies the mathematical analysis.

The linear  model is neutrally stable so the bending and the torsion oscillations do not decay. We seek a compensator to asymptotically stabilize  these oscillations.  A model for designing a stabilizing compensator  does not need to be as accurate  as  a model for simulation because any errors will decay during the stabilization process. 

To the model of \cite{BAH96} we add two  actuators located where the horizontal beam  is joined to its support.  One actuator affects the bending of the beam at its support and other affects the torsion of the beam at its support.   
We first seek a state feedback control law   which drives the oscillations to zero.  To find such a feedback  we set up and solve a Linear Quadratic Regulator (LQR).   We are certainly not the first to use 
LQR to find a state feedback that stabilizes a bending and twisting beam.  Edwards, Breakwell and Bryson \cite{EBB78} used LQR to find a state feedback that stabilizes an airfoil, a cross section of a wing, using  leading edge and trailing edge flaps.   We stabilize the whole beam with two actuators at the root of the beam.

This work was presented at 2025 American Control Conference in Denver, CO and is repeated here because it is necessary background  for the additinal material, the Kalman filter and aerodynamic model, that we present here.
We add two sensors to the model, one which measures the vertical velocity at the tip of the wing and other that measures the rate of torsion also at the tip of the wing.
We design a Kalman filter to process these measurements and return an estimate of the full state.  

Finally we add a model of the aerodynamics forces that  a bending and twisting beam generates to obtain a model of a wing.   The  model is based on the classic Wagner model  for the aerodynamic forces exprienced by an airfoil in flight as described in \cite{BR13}.   We extend a standard finite dimensional  model of an airfoil to a model of the complete wing and expand it using the same families of eigenfunctions that were used to describe the bending and twisting of the beam.

We present  several simulations of the beam and wing to demonstrate the effectiveness of our approach.

\section{Dynamical System}
Let the  $y$ axis be the  axis of rotation of the beam  and suppose it extends from $y=0$ where it is attached to its support and its free end is at $y=L$.  Let $w(y,t) $ be the vertical deflection of beam at location $y$ and time $t$ and let $\theta(y,t)$ be the angle of rotation of the beam around the $y$ axis at  location $y$ and time $t$.   According to \cite{BAH96}, equations (3-155) and (3-156), the free vibrations of a uniform beam are governed by the two inertially coupled linear PDEs
\bea \label{dyn}
 \bmt \mu & -S_y\\ -S_y&  I_y  \emt\bmt \frac{\partial^2 w }{\partial t^2}(y,t)\\ \frac{\partial^2 \theta }{\partial t^2}(y,t)\emt &=&
 \bmt -EI \frac{\partial^4 w }{\partial y^4}(y,t)\\ GJ \frac{\partial^2 \theta }{\partial y^2}(y,t)\emt
\eea
where 
\bean
\begin{array}{llllll}
\mu& \mbox{ mass per unit span}\\
EI& \mbox{  bending rigidity} \\
GJ&\mbox{ torsion rigidity}\\
S_y & \mbox{ static moment per unit span about $y$ axis}\\
I_y & \mbox{ moment of inertia per unit span about $y$ axis}
\end{array} 
\eean
The inertial coupling coefficient $ S_y$ is zero  if all the centers of gravity of the cross sections lie on the elastic axis.

The bending boundary conditions at the free end of the beam are
\bea \label{freebbc}
\frac{\partial^2 w}{\partial y^2}(L,t)=0,&&\frac{\partial^3 w}{\partial y^3}(L,t)=0
\eea
and at the fixed end of the beam we assume that there is an actuator that can deliver a bending moment
 \bea \label{conbbc}
w(0,t)= 0,&& \frac{\partial^2 w}{\partial y^2}(0,t)=B_1u_1(t)
\eea
The torsion boundary condition at the free end of the beam is
\bea \label{freetbc}
\frac{\partial \theta}{\partial y}(L,t)&=& 0
\eea
and at the fixed end of the beam we assume that there is an actuator that can deliver a torque
\bea  \label{contbc}
\frac{\partial \theta}{\partial y}(0,t)&=& B_2u_2(t)
\eea
We define
\bea
B=\bmt B_1&0\\0& B_2\emt
\eea

We wish to express the dynamics as a first order system so we  introduce a four vector valued variable
\bean
z(y,t)&=& \bmt w(y,t)& \frac{\partial w}{\partial t}(y,t)& \theta(y,t) &\frac{\partial \theta}{\partial t}(y,t)\emt\tr 
\eean
then (\ref{dyn}) becomes
\bea \label{dyn1}
M\frac{\partial z}{\partial t}&=&D z(y,t)
\eea
where
\bean
D&=& \bmt 0&1&0&0\\-EI \frac{\partial^4 }{\partial y^4}&0&0&0\\
0&0&0&1\\0&0&GJ  \frac{\partial^2 }{\partial y^2}&0
\emt,
\\
M&=&\bmt 1&0&0&0\\ 0&\mu&0&-S_y\\0&0&1&0\\ 0&-S_y&0&I_y\emt
\eean
\normalsize
The boundary conditions on $z(y,t)$ are
\bea \label{zbc}
z_1(0,t)=0,&& \frac{\partial^2 z_1}{\partial y^2}(0,t) =B_1u_1(t)\\ \nonumber
 \frac{\partial^2  z_1}{\partial y^2}(L,t) =0,&&  \frac{\partial^3  z_1}{\partial y^3}(L,t) =0\\ \nonumber
 \frac{\partial z_3}{\partial y}(0,t) =B_2u_2(t),&&\frac{\partial z_3}{\partial y}(L,t) =0
\eea

We seek a feedback law of the form
\bea
u(t)&=& \int_0^L K(y)z(y,t) \ dy
\eea
to stabilize the bending and torsion oscillations so we set up a Linear Quadratic Regulator (LQR).  We choose 
a   $4\times 4$  nonnegative  definite matrix valued  function $Q(y_1,y_2)$ that is symmetric in its arguments,  $Q(y_1,y_2)=Q(y_2,y_1)$  and a positive definite $ 2\times 2$  matrix $R$.
For a given initial condition $z(y,0)$ we seek to minimize by choice of $u(t)$ the quantity
\bea
\label{crit} \int_0^\infty \iint_{\cal S} z\tr (y_1,t)Q(y_1,y_2)z(y_2,t)+u\tr (t)Ru(t) \ 	dA\ dt
\eea
where ${\cal S}$ is the square $[0,L]^2$ and $dA=dy_1dy_2$.

Let $P(y_1,y_2)$ be an   $4\times 4$  nonnegative  definite matrix valued $C^3$  function that is symmetric in its arguments,  $P(y_1,y_2)=P(y_2,y_1)$ and satisfies these homogeneous boundary conditions for $i,j=1,\ldots,4$ 
\eq \label{Pbc}
P_{2,j}(0,y_2)=0,&&P_{i,2}(y_1,0)=0\\
\frac{\partial^2 P_{2,j}}{\partial  y^2_1}(0,y_2)=0,&&\frac{\partial^2 P_{i,2}}{\partial  y^2_2}(y_1,0)=0\\
\frac{\partial^2 P_{2,j}}{\partial  y^2_1}(L,y_2)=0,&&\frac{\partial^2 P_{i,2}}{\partial  y^2_2}(y_1,L)=0\\
\frac{\partial^3 P_{2,j}}{\partial  y^3_1}(L,y_2)=0,&&\frac{\partial^3 P_{i,2}}{\partial  y^3_2}(y_1,L)=0\\
\frac{\partial P_{4,j}}{\partial  y_1}(0,y_2)=0,&&\frac{\partial P_{i,4}}{\partial  y_2}(y_1,0)=0\\
\frac{\partial P_{4,j}}{\partial  y_1}(L,y_2)=0,&&\frac{\partial P_{i,4}}{\partial  y_2}(y_1,L)=0
\ee
These boundary conditions  are analogous to the those on $z(y,t)$ when the the controls are zero and will be important when we integrate by parts below.  We assume that the $P_{i_1,i_2}(x_1,x_2)$ are $C^3$  so that these boundary conditions make sense.

LQR requires a stabilizability condition so we assume that for each $z(y,0)$ there is a $u(t)$ such that $w(y,t)\to 0$ and $\theta(y,t)\to 0$
 as $t\to \infty$. Then by the Fundamental Theorem of Calculus
\bea
0&=& \iint_{\cal S}z\tr (y_1,0)MP(y_1,y_2)Mz( y_2,0)\ 	dA\\
&&+\int_0^\infty{\partial \over \partial t} \iint_{\cal S}z\tr (y_1,t)MP(y_1,y_2)Mz(y_2,t) \ 	dA \ dt\nonumber
\eea

LQR also requires a detectability condition  so we assume that if $z(y,t)$ is such that 
\bean 
 \int_0^\infty \iint_{\cal S} z\tr (y_1,t)Q(y_1,y_2)z(y_2,t)\ 	dA\ dt&=&0
 \eean
 then $z(y,t)\to 0$ as $t \to \infty$.

We bring the time differentiation inside the spatial integrals and  integrate  by parts several times  to obtain
\tiny
\bea \label{zero}
&&0=  \iint_{\cal S}z\tr (y_1,0)MP(y_1,y_2)Mz(y_2,0) \ 	dA\\
\nonumber&& + \int_0^\infty \iint_{\cal S}z\tr (y_1,t)
\bmt -EI\frac{\partial^4 P_{:,2}}{\partial y_1^4}(y_1,y_2)\\P_{:,1}(y_1,y_2)\\GJ\frac{\partial^2 P_{:,4}}{\partial y_1^2}(y_1,y_2)
\\P_{:,3}(y_1,y_2)
\emt \tr M z(y_2,t)\  dA\ dt \\
&&+ \nonumber\int_0^\infty \iint_{\cal S}z\tr (y_1,t)
M\bmt -EI\frac{\partial^4 P_{:,2}}{\partial y_2^4}(y_1,y_2)\\P_{:,1}(y_1,y_2)\\GJ\frac{\partial^2 P_{:,4}}{\partial y_2^2}(y_1,y_2)
\\P_{:,3}(y_1,y_2)
\emt   z(y_2,t)\  dA\ dt\ \\   \nonumber
&&+\int_0^\infty \int_0^L u\tr (t)B\bmt -EI\frac{\partial P_{:,2}}{\partial y_1}(0,y_2\\-GJP_{:,4}(0,y_2)\emt\tr M z(y_2,t)\ dy_2\ dt\\
&& \nonumber
+\int_0^\infty \int_0^L z\tr (y_1,t) M\bmt -EI\frac{\partial P_{:,2}}{\partial y_2}(y_1,0\\-GJP_{:,4}(y_1,0)\emt Bu(t)\ dy_1 \ dt
\eea
\normalsize
where $P_{i,:}(y_1,y_2)$ and $P_{:,j}(y_1,y_2)$ denote the $i^{th}$ row and the $j^{th}$ column of $P(y_1,y_2)$ respectively.

We add the right side of (\ref{zero}) to the criterion (\ref{crit}) to get an equivalent criterion to be minimized. We
 wish to find a $4 \times 2$ matrix valued function $K(y)$ such that the time integrand of equivalent criterion is equal to a perfect square of the form 
\bean
\iint_{\cal S} \left(u(t)-K(y_1)z(y_1,t)\right)\tr R \left(u(t)-K(y_2)z(y_2,t)\right)\ dA
\eean
The terms quadratic in $u(t)$ match so we equate terms bilinear in $u\tr (t)$ and  $z(y_2,t)$. This yields
\bean
-RK(y_2)&=&B\bmt -EI\frac{\partial P_{2,:}}{\partial y_1}(0,y_2)\\-GJP_{4,:}(0,y_2)\emt M
\eean
so we assume that 
\bean
K(y_2)&=&R^{-1}B\bmt EI\frac{\partial P_{2,:}}{\partial y_1}(0,y_2)\\GJP_{4,:}(0,y_2)\emt M
\eean
By symmetry
\bean
K(y_1)&=&R^{-1}B\bmt EI\frac{\partial P_{2,:}}{\partial y_2}(y_1,0)\\P_{4,:}(y_1.0)\emt M  
\eean
 
 Then by equating terms  bilinear in $z\tr (y_1,t)$ and  $z(y_2,t)$ we obtain the Riccati PDE for quadratic Fredholm kernel $P(y_1,y_2)$  of the optimal cost,
\bea  \nonumber
&&\bmt -EI\frac{\partial^4 P_{2,:}}{\partial y_1^4}(y_1,y_2)\\P_{1,:}(y_1,y_2)\\GJ\frac{\partial^2 P_{4,:}}{\partial y_1^2}(y_1,y_2)\\P_{3,:}(y_1,y_2)\emt \tr M\\
&& \label{Ricpde}
+M\bmt -EI\frac{\partial^4 P_{2,:}}{\partial y_2^4}(y_1,y_2)\\P_{1,:}(y_1,y_2)\\GJ\frac{\partial^2 P_{4,:}}{\partial y_2^2}(y_1,y_2)\\P_{3,:}(y_1,y_2)\emt +Q(y_1,y_2) 
  \\&&=M\bmt EI\frac{\partial P_{2,:}}{\partial y_2}(y_1,0)\\GJP_{4,:}(y_1,0)\emt\tr \Gamma \bmt EI\frac{\partial P_{2,:}}{\partial y_1}(0,y_2)\\
  GJP_{4,:}(0,y_2)\emt \nonumber M
\nonumber
\eea
 \normalsize
where $\Gamma=BR^{-1}B$.
This is an elliptic PDE with a quadratic nonlinearity.  

 \section{Fourier Analysis} 
 We  solve the Riccati PDE (\ref{Ricpde}) by Fourier series.  
 The partial  differential operator $\frac{\partial^2 }{\partial y^2}$
is self adjoint under the zero controlled boundary conditions
\bea \label{wavebc}
\frac{\partial \theta}{\partial y}(0,t)=0,&&
\frac{\partial \theta}{\partial y}(L,t)= 0
\eea
  so all the eigenvalues are real and  the eigenfunctions are orthogonal.   The eigenfunctions are $\Theta_m(y)=\cos {m \pi \over L}y$  and the 	eigenvalues are $\eta_m=-({m\pi\over L})^2$ for $n=0,1,2,\ldots$.
  
 These are not the eigenpairs of the wave equation under the boundary conditions (\ref{wavebc}) but are related to them.  The eigenpairs
  of the wave equation 
  \bean
  \frac{\partial^2 \theta }{\partial t^2}&=& GJ  \frac{\partial^2 \theta }{\partial y^2}
  \eean
 when written as a first order system 
 are $\pm \sqrt{GJ \ }\ {n\pi i\over L}$ which are   strictly imaginary if $m>0$.   The corresponding eigenfunctions are
 \bean
 \bmt \cos {m \pi \over L}y  & \pm \sqrt{GJ\ }\ {m\pi  i\over L}  \cos {m \pi \over L}y
 \emt\tr 
 \eean

The partial differential operator
$ 
-\frac{\partial^4 }{\partial y^4}
$
is self-adjoint when subject to the  boundary conditions 
\bea  \label{bc4}
\phi(0)=0, \frac{\partial^2 \phi}{\partial y^2}(0)=0,
\frac{\partial^2\phi }{\partial y^2}(L)=0, \frac{\partial^3 \phi}{\partial y^3}(L)=0
\eea
Note that
these are not the boundary conditions of a cantilever beam.  For a cantilever beam, the second boundary condition is
$ \frac{\partial \phi}{\partial y}(0)=0$, not $ \frac{\partial^2 \phi}{\partial y^2}(0)=0$

  We look for eigenpairs $ \nu, \Phi(x)$ that satisfy 
 \bea \label{Phi}
 \frac{\partial^4 \Phi}{\partial y^4}(y)&=&  \nu \Phi(y)
\eea
and the boundary conditions (\ref{bc4}).

Also note that   an eigenpair  $ \nu, \Phi(x)$  satisfying (\ref{Phi}) and (\ref{bc4}) is not an eigenpair of beam equation
\bea \label{beam}  
\frac{\partial^2 w }{\partial t^2}(y,t)&=&
-EI \frac{\partial^4 w }{\partial y^4}(y,t)
\eea 
but it is related to them.   When written as a first order system the beam eigenpairs are the eigenpairs of the matrix differential operator
\bean
\bmt 0& 1\\-EI\frac{\partial^4 }{\partial y^4}&0\emt
\eean
The beam eigenvalues are $\pm \sqrt{EI\nu\  \ }$ and since $EI\nu <0$ the beam eigenvalues are strictly imaginary.
The corresponding beam eigenfunctions are vector valued,
\bean 
\bmt \Phi(y)& \pm \sqrt{EI\nu \ \ } \ \Phi(y)
\emt\tr 
\eean
 
Because the differential operator  $-\frac{\partial^4 }{\partial y^4}(x)$
is self-adjoint under the zero controlled boundary conditions (\ref{bc4})  it follows that all its  eigenvalues $\nu_n$  are real and there
is an orthogonal family  of eigenfunctions, $\Phi_n(y)$. Regardless of the boundary conditions, any eigenfunction of  the partial differential operator $-\frac{\partial^4 }{\partial y^4}(x)$ must be of the form 
 \bean \label{ef}
 \Phi(y)&=& a\cos \nu y +b \sin \nu y+c \cosh \nu y + d\sinh \nu y
 \eean
for some real $\nu, a,b,c,d$.  Then the eigenvalue is  $\lambda=-\nu^4$. 
 We can express the boundary equations as a set of homogeneous linear equations depending on $\nu$,
 \tiny
 \bean
 \bmt
 1&0&1&0\\
 -\nu^2 &0& \nu^2&0\\
 -\nu^2 \cos \nu L& -\nu^2 \sin \nu L&\nu^2 \cosh\nu L& \nu^2\sinh \nu L\\
 \nu^3 \sin \nu L& - \nu^3\cos \nu L&  \nu^3\sinh \nu L& \nu^3\cosh \nu L
 \emt
 \bmt a\\b\\c\\d\emt
 = 
 \bmt 0\\0\\0\\0
 \emt
 \eean
 \normalsize
This system has  a nontrivial solution iff $\nu $ is a  root of the determinant of this matrix.	The determinant is
\bean 
-2\nu^7\left(\cos \nu L\sinh \nu L-\cosh \nu L\sin \nu L\right)
\eean

By adjusting the signs of $b$ and $d$ we can restrict our attention to nonnegative roots. 
There is no nonzero eigenfunction corresponding to $\nu=0$ so it is not an eigenvalue.
The positive  roots occur
when $\tan \nu L=\tanh \nu L$.   
   For each $n=1,2,\ldots$ there is exactly one root $ \nu_n    L\in[n\pi,(n+1/2)\pi) $. The first four roots are   $\nu_1 L=3.9266$, $\nu_2 L=7.0686$, $\nu_3 L=10.2102$, $\nu_4 L=13.3518$.  As $n\to\infty$ the $n^{th} $ root  $\nu_n L$ is quickly converging  to $(n\pi +{\pi\over 4})$.  
      
   Because of the boundary conditions
   (\ref{bc4}) at $y=0$ we look for eigenfunctions are of the form
    \bea \label{ef1}
 \Phi(y)&=& b \sin \nu y + d\sinh \nu y
 \eea
   Without loss of generality we can take $b=1$.
Then the boundary conditions
   (\ref{bc4}) at $y=L$ imply
    \bean
  -\sin \nu    L +d \sinh \nu    L&=& 0\\
  - \cos \nu    L +d \cosh \nu    L &=& 0
   \eean
   or equivalently
   \bean
   d&=&{\sin \nu L \over \sinh \nu L}\ =\  {\cos \nu L \over \cosh \nu L}
   \eean
   
   For large $\nu_n    L$ we have $\sinh \nu_n    L\approx\cosh \nu_n    L $
   so we conclude that for large $n$, $ \sin \nu_n    L\approx \cos \nu_n    L $.  
   This happens when $\nu_n    L\approx (n\pi +{\pi\over 4})$. 
 
So the eigenfunctions $\Phi_n (y)$ are converging  to . 
     \bea \label{phim}
  \Phi_n (y)\approx \sin \nu_n   y + {\sin \nu_n L \over \sinh \nu_n L}\sinh \nu y
     \eea
  and 
  \bea \label{phiprimebnd}
  \Phi'_n(0) &=& \nu_n\approx n\pi +{\pi\over 4}
     \eea
     
     We conjecture that 
 the set of vectors  $\{ \Theta_m(y):m=0, 1,2,3,\ldots \}\cup\{\Phi_n(y): n=1,2,3,\ldots \}$ is a Riesz basis for $L^2[0,L]$.

\section{Series Solution of the Riccati PDE}
 Suppose $Q(y_1,y_2)$ 
 has an expansion of the form 
 \small
\bea 
&&Q(y_1,y_2) = \sum_{n=1}^\infty \bmt Q^{n,n}_{1,1} &0&0&0\\0&Q^{n,n}_{2,2}&0&0\\0&0&0&0\\0&0&0&0\emt \Phi_{n}(y_1)\Phi_{n}(y_2) \nonumber
\\&&  +\sum_{m=0}^\infty \bmt 0&0&0&0\\0&0&0&0\\0&0&Q^{m,m}_{3,3} &0\\0&0&0&Q^{m,m}_{4,4}\emt\Theta_m(y_1)\Theta_m(y_2)  \label{Qsum}
\eea

\normalsize
 We could consider more general $Q(y_1,y_2)$ but  to keep the exposition relatively simple we do not.  Notice that the ranges of the sums are different.
 
 We also assume that the solution $P(y_1,y_2)$ of the Riccati PDE (\ref{Ricpde}) has a similar but more complicated expansion.   
 \bean
 P(y_1,y_2)&=& \sum_{n_1,n_2=1}^\infty P^{n_1,n_2}\Phi_{n_1}(y_1)\Phi_{n_2}(y_2)
 \\
 &+&\sum_{n_1=1,m_2=0}^\infty P^{n_1,m_2}\Phi_{n_1}(y_1)\Theta_{m_2}(y_2)
 \\
 &+&\sum_{m_1=0,n_2=1}^\infty P^{m_1,n_2}\Theta_{m_1}(y_1)\Phi_{n_2}(y_2)
 \\
 &+&\sum_{m_1,m_2=0}^\infty P^{m_1,m_2}\Theta_{m_1}(y_1)\Theta_{m_2}(y_2)
 \eean
 Note we are abusing notation, $P^{n_1,n_2}$ is not necessarily the same as $P^{m_1,m_2}$
 even when $n_1=m_1$ and $n_2=m_2$. We use $n$ as a superscipt to indicate a coefficient
 of $\Phi_n(y)$ and we use $m$ as a superscipt to indicate a coefficient
 of $\Theta_m(y)$.

We plug these expansions into Riccati PDE 
and collect similar terms to obtain an infinite dimensional algebraic Riccati equation which has four uncoupled components.
 The
 $ \Phi_{n_1}(y_1)\Phi_{n_2}(y_2)$ component is 
 \bea \label{arenn}  \nonumber
 &&\bmt \lambda_{n_1} EIP^{n_1,n_2}_{2,:}\\P^{n_1,n_2}_{1,:}\\ -\eta_{n_1}^2GJP^{n_1,n_2}_{4,:}\\P^{n_1,n_2}_{3,:}\emt\tr M
 +M\bmt \lambda_{n_2} EIP^{n_1,n_2}_{2,:}\\P^{n_1,n_2}_{1,:}\\ -\eta_{n_2}^2GJP^{n_1,n_2}_{4,:}\\P^{n_1,n_2}_{3,:}\emt\\
\label{arenn}
&&+ \delta_{n_1,n_2}\bmt Q^{n_1,n_1}_{1,1}&0&0&0\\0&Q^{n_1,n_1}_{2,2}&0&0&\\0&0&0&0\emt\\
&&=M\bmt EI \sum_{n_4=1}^\infty P^{n_1,n_4}_{2,:} \Phi'_{n_4}(0)\\  GJ \sum_{n_4=1}^\infty P^{n_1,n_4}_{4,:} \Phi_{n_4}(0) \emt\tr \Gamma \nonumber\\ \nonumber
&&\times  \bmt EI \sum_{n_3=1}^\infty P^{n_3,n_2}_{2,:} \Phi'_{n_3}(0)\\  GJ \sum_{n_3=1}^\infty P^{n_3,n_2}_{4,:} \Phi_{n_3}(0) \emt M
 \eea
 
  The
$ \Phi_{n_1}(y_1)\Theta_{m_2}(y_2)$  component is 
\bea \nonumber 
 &&\bmt \lambda_{n_1} EIP^{n_1,n_2}_{2,:}\\P^{n_1,n_2}_{1,:}\\ -\eta_{n_1}^2GJP^{n_1,n_2}_{4,:}\\P^{n_1,n_2}_{3,:}\emt\tr  M+
M\bmt  -\nu^2_{m_2} EI P^{m_1,m_2}_{2,:} \\ P^{m_1,m_2}_{1,:} \\   \nu_{m_2}GJ P^{m_1,m_2}_{4,:} \\ P^{m_1,m_2}_{3,:} \emt\\
&& \nonumber=M\bmt\sum_{m_4=0}^\infty EI P^{n_1,m_4}\Theta'_{m_4}(0)\\ \sum_{m_4=0}^\infty GJP_{4,:}^{n_1,m_4}\Theta_{m_4}(0)\emt\tr\Gamma \\
&&\times \bmt \sum_{n_3=1}^\infty EI P_{2,:}^{n_3,m_2} \Phi'_{n_3}(0)  \\
 \sum_{n_3=1}^\infty GJ P^{n_1,m_4}_{4,:} \Phi_{n_3}(0)\emt M\label{arenm}
\eea
A simple solution to this equation is to take $P^{n,m}=0$ for all
$n$ and $m$.  Then by symmetry $P^{m,n}=0$.  This would not be possible if we had chosen $Q^{m,n}\ne 0$.

Finally  the
$ \Theta_{m_1}(y_1)\Theta_{m_2}(y_2)$  component is 
\bea \nonumber 
&&\bmt  -\nu^2_{m_1} EI P^{m_1,m_2}_{2,:} \\ P^{m_1,m_2}_{1,:} \\   \nu_{m_1}GJ P^{m_1,m_2}_{4,:} \\ P^{m_1,m_2}_{3,:} \emt\tr M +
M\bmt  -\nu^2_{m_2} EI P^{m_1,m_2}_{2,:} \\ P^{m_1,m_2}_{1,:} \\   \nu_{m_2}GJ P^{m_1,m_2}_{4,:} \\ P^{m_1,m_2}_{3,:} \emt\\
&&+ \delta_{m_1,m_2}\bmt 0&0&0&0\\0&0&0&0\\0&0&Q^{m_1,m_1}_{3,3}&0\\0&0&0&Q^{m_1,m_1}_{4,4}\emt\label{aremm} \\
&&=M\bmt 0^{1\times 4}\\  GJ \sum_{m_4=1}^\infty P^{m_1,m_4}_{4,:}  \emt\tr \Gamma\nonumber\\&& \nonumber
\times \bmt 0^{1\times 4}\\  GJ \sum_{n_3=1}^\infty P^{m_3,m_2}_{4,:}  \emt M
 \eea

 \section{Policy Iteration}
   We can approximately solve the algebraic Riccati equations  (\ref{arenn}) and  (\ref{aremm}) by policy iteration.  The method can also be seen as value iteration. 
  To find an initial iterate $\left(P^{n_1,n_2}_{i,j}\right)^{(0)}$ and $ \left(P^{m_1,m_2}_{i,j}\right)^{(0)}$ we make the following assumptions.
  \begin{enumerate}
  \item  $\left(P^{n_1,n_2}_{i,j}\right)^{(0)}=0$ if $i=3,4$ or  $j=3,4$.
  \item $\left(P^{n_1,n_2}_{i,j}\right)^{(0)}=0$ if $n_1\ne n_2$ when $i,j=1,2$.
  \item  $\left(P^{m_1,m_2}_{i,j}\right)^{(0)}=0$ if $i=1,2$ or  $j=1,2$.
  \item $\left(P^{m_1,m_2}_{i,j}\right)^{(0)}=0$ if $m_1\ne m_2$ when $i,j=3,4$.
 \end{enumerate}

 These assumptions decouple
  bending from twisting so the four dimensional equations (\ref{arenn}) and (\ref{aremm}) reduce to two independent  two dimensional equations.
The Fourier expansion of  the bending model  is solely in terms $\Phi_n(y)$ and it satisfies     
\tiny
\bea \label{arenns}
&& \bmt \lambda_{n_1} \mu EIP^{n_1,n_2}_{1,2}&P^{n_1,n_2}_{1,1}\\\lambda_{n_1} \mu EIP^{n_1,n_2}_{2,2}&P^{n_1,n_2}_{1,2}\emt\\
&&  \nonumber
 +\bmt \lambda_{n_2} \mu EI P^{n_1,n_2}_{2,1}&\lambda_{n_2} \mu EI P^{n_1,n_2}_{2,2}\\P^{n_1,n_2}_{1,1}&P^{n_1,n_2}_{1,2}\emt
 \nonumber\\  \nonumber
 &&
+ \delta_{n_1,n_2}\bmt Q^{n_1,n_1}_{1,1}&0\\0&Q^{n_1,n_1}_{2,2}\emt\\
&&=\bmt \sum_{n_4=1}^\infty \mu EI  P^{n_1,n_4}_{1,2} \Phi'_{n_4}(0)\\\sum_{n_4=1}^\infty EI  P^{n_1,n_4}_{2,2} \Phi'_{n_4}(0)\emt \Gamma_{1,1} \nonumber\\ \nonumber
&&\times  \bmt \sum_{n_3=1}^\infty \mu EI  P^{n_3,n_2}_{2,1} \Phi'_{n_3}(0)&\sum_{n_3=1}^\infty  EI P_{2,2}^{n_3,n_2}\Phi'_{n_3}(0)\emt
 \eea
 \normalsize
 Recall $\Phi'_{n}(0)\to \nu_n$ as $n\to \infty$.
 
 The Fourier expansion of  the twisting model  is solely in terms of $\Theta_m(y)$ and it satisfies
   \bea \label{aremms}
 &&  \bmt \eta_{m_1}GJ P_{3,4}^{m_1,m_2}& I_y P_{3,3}^{m_1,m_2}\\ \eta_{m_1}GJ P_{4,4}^{m_1,m_2}& I_y P_{3,4}^{m_1,m_2}\emt
   \\ \nonumber
   && + \bmt \eta_{m_2}GJ P_{4,3}^{m_1,m_2}&\eta_{m_2}GJ P_{4,4}^{m_1,m_2}\\I_yP^{m_1,m_2}_{3,3}&I_yP^{m_1,m_2}_{3,4}\emt\\  \nonumber
   && +\delta_{m_1,m_2}\bmt Q^{m_1,m_1}_{3,3}&0\\ 0& Q^{m_1,m_1}_{4,4}  \emt
   \\  \nonumber 
   && =\bmt GJ \sum_{m_4=0}^\infty P^{m_1,m_4}_{3,4}\\ I_y GJ \sum_{m_4=0}^\infty P^{m_1,m_4}_{4,4}\emt \Gamma_{2,2}
  \\
  &&\times
   \bmt GJ\sum_{m_3=0}^\infty P^{m_1,m_2}_{4,3}& I_y GJ \sum_{m_3=0}^\infty P^{m_3,m_2}_{4,4}\emt \nonumber
 \eea

 The above assumptions  on the initial iterate simplifies (\ref{arenns})   to
 \bea \label{arennss}
&& \bmt \lambda_{n} \mu EI\left(P^{n,n}_{1,2}\right)^{(0)}&\left(P^{n,n}_{1,1}\right)^{(0)}\\\lambda_{n} \mu EI\left(P^{n,n}_{2,2}\right)^{(0)}&\left(P^{n,n}_{1,2}\right)^{(0)}\emt \nonumber
\\&& \nonumber+\bmt \lambda_{n } \mu EI \left(P^{n,n}_{2,1}\right)^{(0)}&\lambda_{n} \mu EI \left(P^{n,n}_{2,2}\right)^{(0)}\\\left(P^{n,n}_{1,1}\right)^{(0)}&\left(P^{n,n}_{1,2}\right)^{(0)}\emt
 \nonumber\\ 
 &&
+ \bmt Q^{n,n}_{1,1}&0\\0&Q^{n,n}_{2,2}\emt\\
&&=\bmt \mu EI  \left(P^{n,n}_{1,2}\right)^{(0)} \Phi'_{n}(0)\\ EI  \left(P^{n,n}_{2,2}\right)^{(0)} \Phi'_{n}(0)\emt \Gamma_{1,1} \nonumber\\ \nonumber
&&\times  \bmt \mu EI  \left(P^{n,n}_{2,1}\right)^{(0)} \Phi'_{n}(0)&  EI \left(P_{2,2}^{n,n}\right)^{(0)}\Phi'_{n}(0)\emt
 \eea
 and  simplifies   (\ref{aremms}) to
 \bea \label{aremmss}
 &&  \bmt \eta_{m}GJ \left(P_{3,4}^{m,m}\right)^{(0)}& I_y \left(P_{3,3}^{m,m}\right)^{(0)}\\ \eta_{m}GJ \left(P_{4,4}^{m,m}\right)^{(0)}& I_y \left(P_{3,4}^{m,m}\right)^{(0)}\emt
   \\ \nonumber
   && + \bmt \eta_{m }GJ \left(P_{4,3}^{m,m}\right)^{(0)}&\eta_{m}GJ \left(P_{4,4}^{m,m}\right)^{(0)}\\I_y\left(P^{m,m}_{3,3}\right)^{(0)}&I_y\left(P^{m,m}_{3,4}\right)^{(0)}\emt\\  \nonumber
   && +\bmt Q^{m,m}_{3,3}&0\\ 0& Q^{m,m}_{4,4}  \emt
   \\  \nonumber 
   && =\bmt GJ \left(P^{m,m}_{3,4}\right)^{(0)}\\ I_y GJ  \left(P^{m,m}_{4,4}\right)^{(0)}\emt \Gamma_{2,2}
  \\
  &&\times
   \bmt GJ \left(P^{m,m}_{4,3}\right)^{(0)}& I_y GJ  \left(P^{m,m}_{4,4}\right)^{(0)}\emt \nonumber
 \eea
    
  Notice the $1,1$ component of (\ref{arennss}) is a quadratic equation in $\left(P_{1,2}^{n,n}\right)^{(0)}=\left(P_{2,1}^{n,n}\right)^{(0)}$,
  \bea
  2\lambda_n EI \left(P_{1,2}^{n,n}\right)^{(0)}+Q^{n,n}_{1,1}=\Gamma_{1,1}\left(EI \left(P^{n,n}_{1,2}\right)^{(0)}\Phi'_{n}(0)\right)^2
  \eea
  Since it is harder to stabilize when $z_1(y,0)$ and $z_2(y,0)$ have the same sign we take the positive
  root of this quadratic to be the initial estimate 
  \bea \label{P120}
  \left(P_{1,2}^{n,n}\right)^{(0)}={\lambda_nEI+\sqrt{(\lambda_nEI)^2+Q^{n,n}_{1,1}\Gamma_{1,1}\left(EI\Phi'_{n}(0)\right)^2}
  \over \Gamma_{1,1}\left(EI\Phi'_{n}(0)\right)^2}
  \eea
     
  The $2,2$ component of (\ref{arennss}) is a quadratic equation in $\left(P^{n,n}_{2,2}\right)^{(0)}$
  \bea
  2I_y \left(P_{1,2}^{n,n} \right)^{(0)}+Q^{n,n}_{2,2}=\Gamma_{2,2}\left(I_y EI \left(P^{n,n}_{2,2}\right)^{(0)} \Phi'_n(0)\right)^2
  \eea 
  and since $\left(P^{n,n}_{2,2}\right)^{(0)}$ must be nonnegative we define 
  \bea
  \label{P220} 
   \left(P_{2,2}^{n,n}\right)^{(0)}&=&{\sqrt{2I_y \left(P_{1,2}^{n,n}\right)^{(0)} +Q^{n,n}_{2,2}}\over \Gamma_{2,2}^{1/2}I_y EI  \Phi'_n(0)}
  \eea
  
  We can use the $1,2$ or the $2,1$ component of (\ref{arennss}) to define $ \left(P_{1,1}^{n,n}\right)^{(0)}$.  Both lead to
  \bea \label{P110}
  &&\left(P_{1,1}^{n,n}\right)^{(0)}= -\lambda_n\mu I_yEI \left(P_{2,2}^{n,n}\right)^{(0)}\\
  &&+\mu EI \left(P_{1,2}^{n,n}\right)^{(0)}\Gamma_{1,1}I_yEI   \left(P_{2,2}^{n,n}\right)^{(0)}\Phi'_n(0)\nonumber
  \eea
  
  Then our initial iterate of the kernel of the  bending part of the  optimal  cost is 
  \bea
  \label{Psumnn}
   \sum_{n=1}^\infty \bmt \left(P_{1,1}^{n,n}\right)^{(0)}&\left(P_{1,2}^{n,n}\right)^{(0)}\\
  \left(P_{2,1}^{n,n}\right)^{(0)}&\left(P_{2,2}^{n,n}\right)^{(0)}\emt \Phi_n(y_1)\Phi_n(y_2) 
   \eea
  
  \noindent
 {\bf Theorem 1:} 
 The series (\ref{Psumnn}) converges to a continuous $2\times 2$ matrix valued function.   if there 
 exist  numbers $q>0$ and $r>8$ such that 
 \bea \label{q1r1}
 \left|Q_{i,j}\right|&\le& {q\over m^{r}}
 \eea 
 or $i,j=1,2$.\\
 In addition if for some  integer $\rho\ge 0$ 
 \bea \label{q1r1rho}
 \left|Q_{i,j}\right|&\le& {q\over m^{r+\rho}}
 \eea 
 or $i,j=1,2$. then the series (\ref{Psumnn}) converges to a $C^{\rho}$  $2\times 2$ matrix valued function.\\
 {\bf Proof:}
 The Mean Value Theorem applied to the function $g(s)=\sqrt{s}$ on the interval $[a,b]$ implies that there
 is an $s$ in $(a,b)$ such that 
 \bean
 \sqrt{b}-\sqrt{a}={1\over 2\sqrt{s}} (b-a)
 \eean
 Let $a=(\lambda_nEI)^2$ and $b=(\lambda_nEI)^2 +Q_{1,1}(EI\Phi'_n (0))^2\Gamma_{1,1}$. Then
 $ \lambda_nEI=-\sqrt{(\lambda_nEI)^2}$ 
 so  there exists an $s$ between $(\lambda_nEI)^2$ and 
 $(\lambda_nEI)^2 +Q_{1,1}(EI\Phi'_n (0))^2\Gamma_{1,1}$ such that 
 \bean
 \left(P_{1,2}^{n,n}\right)^{(0)}&=& {1\over 2\sqrt{s}}Q_{1,1}^{n,n}
 \eean
 The maximum of the right side of this equation  between $(\lambda_nEI)^2$ and 
 $(\lambda_nEI)^2 +Q_{1,1}(EI\Phi_n (0)^2\Gamma_{1,1}$  occurs at 
 $s=(\lambda_nEI)^2$
 so 
 \bean
 \left|\left(P_{1,2}^{n,n}\right)^{(0)}\right|&\le & {Q_{1,1}^{n,n}\over 2\left|\lambda_n   \right|EI}
 \eean

Since $\left|\lambda_n   \right| =\nu_n^4\approx \left({ n\pi +\pi/4\over L}\right)^4$,  $Q^{n,n}_{1,1}\le {q\over n^r}$  with $r>8$ and $\left|\Phi_n (y)\right|\le 2$ the series 
 \bean
  \sum_{m=0}^\infty \left(P_{1,2}^{n,n}\right)^{(0)}\Phi_n (y_1)\Phi_n (y_2)
 \eean 
 converges uniformly to a continuous function $P^{(0)}_{1,2}(y_1,y_2)$.

  Now by (\ref{phiprimebnd}) $\Phi_n '(0)=\nu_n    = O(n)$ so
 \bean
 \left(P_{2,2}^{n,n}\right)^{(0)}&=& {\sqrt{2I_y\left(P_{1,2}^{n,n}\right)^{(0)}+Q_{2,2}^{n,n}}\over \Gamma_{2,2}^{1/2} I_yEI\Phi'_n(0)}
 \eean
which  is $O(n^{-1-r/2})$.   Since $r>8$ the series 
  \bean
  \sum_{n=0}^\infty \left(P_{2,2}^{n,n}\right)^{(0)}\Phi_n (y_1)\Phi_n (y_2)
 \eean 
 converges uniformly to a continuous function $P^{(0)}_{2,2}(y_1,y_2)$.

Finally 
\bean
 &&\left(P_{1,1}^{n,n}\right)^{(0)}=-\lambda_n    I_y EI \left(P_{2,2}^{n,n}\right)^{(0)}\\
 &&+\left(EI\Phi'_{n}(0)\right)^2I_y\left(P_{1,2}^{n,n}\right)^{(0)}\left(P_{2,2}^{n,n}\right)^{(0)}
 \eean
 so $\left(P_{1,1}^{n,n}\right)^{(0)}=O(n^{3-r/2})$ and the
 series 
  \bean
  \sum_{m=0}^\infty \left(P_{1,1}^{n,n}\right)^{(0)}\Phi_n (y_1)\Phi_n (y_2)
 \eean 
 converges uniformly to a continuous function $P^{(0)}_{1,1}(y_1,y_2)$ because  $3-r/2<-1$.\\
 \noindent
 The assertion (\ref{q1r1rho}) is obtained by term by term differentiation of the series (\ref{Psumnn}).
{\bf QED}

   The way we interpret this result is that if we want the  cost  of moving the infinite number of imaginary eigenvalues  into the open left half plane to be finite we must make $Q_{i,i}^{n,n}$ go to zero quite rapidly.
   
   Notice the $3,3$ component of (\ref{aremmss}) is a quadratic equation in $P_{3,4}^{n,n}=P_{4,3}^{n,n}$
  \bea
  2\eta_mGJ P_{3,4}^{m,m}+Q^{m,m}_{3,3}=\Gamma_{2,2}\left(GJ P^{m,m}_{3,4}\right)^2
  \eea
  Since it is harder to stabilize when $z_3(y,0)$ and $z_4(y,0)$ have the same sign we take the positive
  root of this quadratic to be the initial estimate 
  \bea \label{P340}
  \left(P_{3,4}^{m,m}\right)^{(0)}={\eta_mGJ+\sqrt{(\eta_mGJ)^2+Q^{m,m}_{3,3}\Gamma_{2,2}\left(GJ\right)^2}
  \over \Gamma_{2,2}\left(GJ\right)^2}
  \eea
     
  The $4,4$ component of (\ref{aremmss}) is a quadratic equation in $P^{m,m}_{4,4}$
  \bea
  2I_y P_{3,4}^{m,m} +Q^{m,m}_{4,4}&=&\Gamma_{2,2}\left(I_y GJ P^{m,m}_{4,4}\right)^2
  \eea 
  and since $P^{m,m}_{4,4}$ must be nonnegative we define 
  \bea
  \label{P440} 
   \left(P_{4,4}^{m,m}\right)^{(0)}&=&{\sqrt{2I_y \left(P_{3,4}^{m,m}\right)^{(0)} +Q^{m,m}_{4,4}}\over \Gamma_{2,2}^{1/2}I_y GJ }
  \eea
  
  We can use the $3,4$ or the $4,3$ component of (\ref{aremmss}) to define $ \left(P_{3,3}^{m,m}\right)^{(0)}$.  Both lead to
  \bea \label{P330}
  &&\left(P_{3,3}^{m,m}\right)^{(0)}= { 1\over I_y}\left(\eta_m GJ \left(P_{4,4}^{m,m}\right)^{(0)}\right.\\
  && \left.+(GJ)^2\left(P_{3,4}^{m,m}\right)^{(0)}\Gamma_{1,1}  \left(P_{4,4}^{m,m}\right)^{(0)}\right)\nonumber
  \eea

Our initial iterate of the kernel of the  twisting part of the  optimal  cost is 
  \bea
  \label{Psummm}
 \sum_{m=0}^\infty \bmt \left(P_{3,3}^{m,m}\right)^{(0)}&\left(P_{3,4}^{m,m}\right)^{(0)}\\
  \left(P_{4,3}^{m,m}\right)^{(0)}&\left(P_{4,4}^{m,m}\right)^{(0)}\emt \Theta_m(y_1)\Theta_m(y_2) 
  \eea
  
  \noindent
 {\bf Theorem  2:} 
 The series (\ref{Psummm}) converges to a continuous $2\times 2$ matrix valued function   if there 
 exist  numbers $q>0$ and $r>6$ such that for $m>0$
 \bea \label{q1r2}
 \left|Q_{i,j}\right|&\le& {q\over m^{r}}
 \eea 
  $i,j=3,4$.\\
   In addition if for some integer $\rho\ge 0$ 
 \bea \label{q1r2rho}
 \left|Q_{i,j}\right|&\le& {q\over m^{r+\rho}}
 \eea 
 or $i,j=1,2$. then the series (\ref{Psummm}) converges to a $C^{\rho}$  $2\times 2$ matrix valued function.\\
 {\bf Proof:}
 Since $\eta_m = -\left({m\pi\over L}\right)^2$ the Mean Value Theorem applied to (\ref{P340}) implies that there is an $s$ 
 between $\left(\eta_m GJ\right)^2$ and $\left(\eta_m GJ\right)^2+Q_{3,3}^{m,m}\Gamma_{2,2}\left(GJ\right)^2$ such that 
 \bea 
 \left(P_{3,4}^{m,m}\right)^{(0)}&=& {1\over 2\sqrt{s}}Q_{3,3}^{m,m}
 \eea
 The maximum of the right side of this equation occurs at $s=\left(\eta_m GJ\right)^2$ so we conclude that
 \bea
 0\le   \left(P_{3,4}^{m,m}\right)^{(0)} \le {1\over 2\left|\eta_m GJ\right|}Q_{3,3}^{m,m}
 \eea
and $\left(P_{3,4}^{m,m}\right)^{(0)}$ is of order ${ 1\over m^8}$.  Since $\left|\Theta_m(y)\right|\le 1$ the series 
 \bea
 P^{(0)}_{3,4}(y_1,y_2)&=& \sum_{m=0}^\infty \left(P_{3,4}^{m,m}\right)^{(0)} \Theta_m(y_1)\Theta_m(y_2)
 \eea
 converges uniformly to a continuous function. 
 
 Now from (\ref{P440}) we know that $\left(P_{4,4}^{m,m}\right)^{(0)}$ is of order ${ 1\over m^{r/2}}$ so 
 the series 
 \bea
 P^{(0)}_{4,4}(y_1,y_2)&=& \sum_{m=0}^\infty \left(P_{4,4}^{m,m}\right)^{(0)} \Theta_m(y_1)\Theta_m(y_2)
 \eea
 converges uniformly to a continuous function. 
 
 Finally from (\ref{P330}) we see that $\left(P_{4,4}^{m,m}\right)^{(0)}$ is of order ${ 1\over m^{r/2-2}}$
and since $r>6$ we have ${r\over 2}-2 >1$ so 
 the series 
 \bea
 P^{(0)}_{3,3}(y_1,y_2)&=& \sum_{m=0}^\infty \left(P_{3,3}^{m,m}\right)^{(0)} \Theta_m(y_1)\Theta_m(y_2)
 \eea
 converges uniformly to a continuous function. \\

  The statement (\ref{q1r2rho}) is obtained by term by term differentiation of the series (\ref{Psummm}).
\noindent
{\bf QED}

 Succesive iterates  are found by plugging $ \left(P^{n_1,n_2}_{i,j}\right)^{(k)}$  and 
$ \left(P^{m_1,m_2}_{i,j}\right)^{(k)}$ into the right sides
of the algebraic Riccati equations (\ref{arenn}) and  (\ref{aremm}) and  solving for $ \left(P^{n_1,n_2}_{i,j}\right)^{(k+1)}$ 
and $ \left(P^{m_1,m_2}_{i,j}\right)^{(k+1)}$ on  the left  side.   A complication arises when $m_1=m_2=0$ because 
then $ \nu_0=0$ and so the left sides is not an invertible expression in $ \left(P^{0,0}_{2,:}\right)^{(k+1)}$
and $ \left(P^{0,0}_{4,:}\right)^{(k+1)}$.    A way around this is for all $m_1^2+m_2^2>0$ to plug $ \left(P^{n_1,n_2}_{i,j}\right)^{(k)}$  and 
$ \left(P^{m_1,m_2}_{i,j}\right)^{(k)}$  into the right sides
of the algebraic Riccati equations (\ref{arenn}) and  (\ref{aremm}) and  solve for  $ \left(P^{n_1,n_2}_{i,j}\right)^{(k+1)}$,  $ \left(P^{m_1,m_2}_{2,:}\right)^{(k+1)}$
and $ \left(P^{m_1,m_2}_{4,:}\right)^{(k+1)}$ for all $n_1,n_2$ and all $m_1,m_2$ such that $m_1^2+m_2^2>0$.  Then plug these $k+1$ iterates into (\ref{arenn}) and  (\ref{aremm}) and solve the resulting quadratics for $ \left(P^{m,m}_{2,:}\right)^{(k+1)}$
and $ \left(P^{m,m}_{4,:}\right)^{(k+1)}$ when $m_1=m_2=0$.  If the resulting quadratics do not determine  $ \left(P^{m,m}_{2,:}\right)^{(k+1)}$
and $ \left(P^{m,m}_{4,:}\right)^{(k+1)}$ when $m_1=m_2=0$ just set them to zero.
It should be noted that if $S_y\ne 0$ then the $k^{th}$ approximations to the bending and twisting kernels are not $2\times 2$ like 
(\ref{Psumnn}) and (\ref{Psummm}) but are
 $4\times 4$ when $k>0$.

This is value (or policy) iteration so we know that the iterates are nonincreasing. If hypotheses of Theorems 1 and 2  are satisfied then the initial approximations (\ref{Psumnn}) and (\ref{Psummm}) are continuous hence bounded on $[0,L]\times [0,L]$.   Hence the initial estimate of the optimal
cost is bounded as are all successive estimates.   From this we conclude that all of the approximations to the optimal feedback move all the closed loop eigenvalues into the open left half plane.
  But we cannot conclude that the closed loop system is exponentially stable because the closed loop eigenvalues could (and do) approach the imaginary axis as the wave numbers $n$ and $m$ increase.
 
 \section{Finite Dimensional Approximating LQR}
 We construct a finite dimensional LQR whose algebraic  Riccati equation is a truncation of (\ref{arenn}) and (\ref{aremm}).
 We choose an $N>0$  and construct a linear system  with state $\zeta=[\zeta_1,\zeta_2,\zeta_3,\zeta_4]\tr $ where 
 $\zeta_1,\zeta_2,\zeta_3,\zeta_3$ is each of dimension $N\times 1$.
  So the finite dimensional state $ \zeta$ is of dimension $4N\times 1$.
  The dynamics is 
 \bean 
\dot{ \zeta}_1&=& \zeta_2 \\
m\dot{ \zeta}_2-S_y\dot{\zeta}_4 &=& F_1\zeta_1 +G_1u_1\\
\dot{ \zeta}_3&=& \zeta_4 \\
-S_y\dot{ \zeta}_2+I_y\dot{\zeta}_4 &=& F_2\zeta_3 +G_2u_2\\
 \eean
 where
 \bean
 F_1=EI\bmt \nu_1&\ldots &0\\ 
 & \ddots& \\
 0&\ldots& \nu_N
  \emt,&& G_1  =EI\bmt 
   \Phi'_{1}(0)\\ \vdots \\ \Phi'_N(0)\emt
  \\
  F_2=GJ\bmt \eta_0&\ldots &0\\ 
 & \ddots& \\
 0&\ldots& \eta_{N-1}
  \emt,&& G_2  =GJ \bmt  \Theta_0(0)\\ \vdots \\  \Theta_{N-1}(0)\emt
 \eean

 This finite dimensional system approximates the infinite dimensional  system in the following manner
\bean
 z_1(y,t)&\approx &\bmt \Phi_1(y)&\ldots & \Phi_N(y)\emt \zeta_1(t)\\
  z_2(y,t) &\approx &\bmt \Phi_1(y)&\ldots & \Phi_N(y)\emt \zeta_2(t)\\
   z_3(y,t)&\approx &\bmt \Theta_0(y)&\ldots & \Theta_{N-1}(y)\emt \zeta_3(t)\\
  z_4(y,t) &\approx &\bmt \Theta_0(y)&\ldots & \Theta_{N-1}(y)\emt \zeta_4(t)
\eean
 
 \section{Example}
 We consider a $N=4$  example which leads to a $16$ dimensional approximation.   For the time being we take all constants equal $1$ except $S_y=1/2$, $Q$ a $16\times 16$ identity matrix and $R$ a $2\times 2$ identity matrix.  The $16$ open and closed loop poles are  
 \bean
 \bmt
 \mbox{ Open Loop Poles}& \mbox{Closed Loop Poles}\\
\pm0.0361i&-0.0114\pm0.0350i\\
\pm0.2067i&-0.0165\pm 0.2051i\\
\pm 0.0724i&-0.0185\pm 0.0724i\\
\pm0.1087ii&-0.0278\pm 0.1087i\\
\pm 0.1450i&-0.0296\pm0.1325i\\
\pm0.0665i&-0.0429\pm0.6672i\\
\pm1.3914i&-0.0652\pm1.3911i\\
\pm 2.3784&-0.0874\pm 2.3700i
\emt \times 10^2
 \eean
 Notice how close the imaginary parts of the closed loop poles are  to the open loop poles, All the closed loop poles have real parts less than $-1$.
 
\begin{figure}
 \centering
\includegraphics[width=4in]{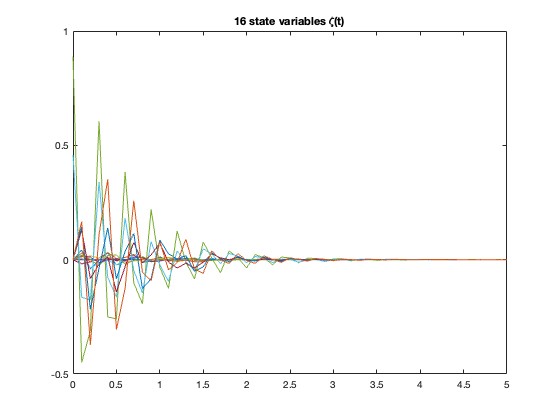}
\caption{Finite Dimensional State $\zeta(t)$  Under Full State Feedback}
\label{fds}        
 \end{figure}

Figure \ref{fds} is a  simulation of the $16$ dimesional state under full state feedabck
    
\begin{figure}
 \centering
\includegraphics[width=3.5in]{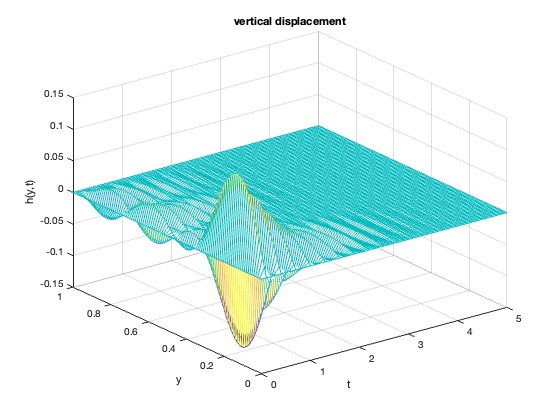}
\caption{Vertical Displacement $h(y,t)$  Under Full State Feedback}    
\label{vdfsf}    
 \end{figure}

    Figure \ref{vdfsf}  shows the vertical displacement
    converging to zero under full state feedback.    
The control input at the root of the beam causes the ripples while stabilizing the vertical displacement of the  beam.

\begin{figure}
 \centering
\includegraphics[width=3.5in]{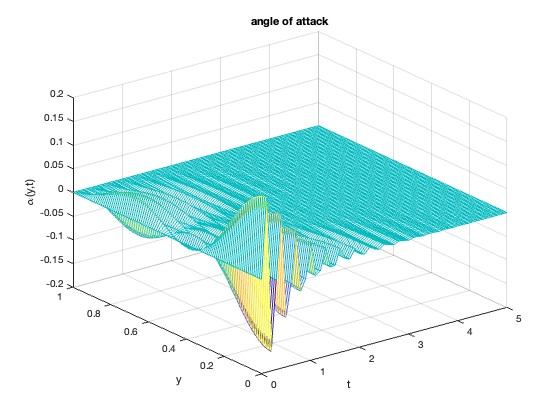}
\caption{Angle of Attack $\alpha(y,t)$ Under Full State Feedback}    
\label{aoafsf}    
 \end{figure}

Figure \ref{aoafsf} shows the stabilization of the angle of attack under full state feedback.
 Again the control input at the root of the beam causes the ripples while stabilizing the torsion of the  beam.

  \section{Kalman Filter}
 In the above sections we have constructed a full state feedback to asymptotically stabilize the linear system.   But we don't have the ability to measure the full state so need to construct a Kalman filter to deliver an estimate of the full state from a finite number of point measurements.    We assume that we can measure 
 $z_2(L,t) =\frac{\partial h}{\partial t}(L,t)$ and $z_4(L,t) =\frac{\partial \alpha}{\partial t}(L,t)$. We could use more sensors and/or different types of sensors but to keep thing simple we do not do so.
 
 To develop the Kalman filter it is simpler to express the inertially coupled bending and twisting dynamics in a more standard form by multiplying the dynamics by the inverse of the coupling matrix $M$.  We also add a distributed driving noise input $v(t) =[v_1(t), v_2(t)]\tr $ and obtain
 \bea \label{dyn_n}
  \frac{\partial z}{\partial t}(y,t)&=& {\cal A} z(y,t) +{\cal B}(y)v(t)
 \eea
 where the matrix partial differential operator is ${\cal A}=M^{-1}D$ and 
 \bean
 {\cal B}(y)&=& \bmt 0&0\\{\cal B}_{11}(y)&{\cal B}_{12}(y)\\ 0&0\\{\cal B}_{21}(y)&{\cal B}_{22}(y)\emt
 \eean
 The boundary conditions on $z(y,t)$ are
  \bea \label{zzbc}
 z_1(0,t)=0, && \frac{\partial^2 z_1}{\partial  y^2}(0,t)=B_1u_(t)\\  \nonumber
 \frac{\partial^2 z_1}{\partial  y^2}(L,t)=0,&& \frac{\partial^3 z_1}{\partial  y^3}(L,t)=0\\  \nonumber
  \frac{\partial z_3}{\partial  y}(0,t)=B_2u_2(t),&& \frac{\partial z_3}{\partial  y}(L,t)=0
 \eea

 The measurements are the vertical velocity $\psi_1(t)=z_2(L,t)$ and angular velocity $\psi_2(t)=z_4(L,t)$ of the tip of the beam corrupted by white Gaussian noise,
    \bean
     \psi(t)  &=&\int_0^L{\cal C}(y) z(y,t)\ dy+ {\cal D}\bmt  w_1(t)\\ w_2(t)\emt
   \eean
  where 
  \bean
  {\cal C}(y)&=&\bmt 0 &\delta(y-L)& 0&0\\0&0&0 &\delta(y-L)\emt
  \eean
  and ${\cal D}$ is a $2\times 2$ invertible matrix  with $w_1(t), w_2(t) $   standard independent white Gaussian noises.   One way of obtaining such measurements is to integrate a pair of accelerometers as was done by Banks et al.\cite{BS97}.
 
  Because of the Gaussian and linear assumptions we expect that when $u_1(t)=u_2(t)=0$ the optimal least squares estimate  $\hat{z}(y,t)$ of $z(y,t) $ is a linear functional of the past observations,
  \bean
 &&\hat{z}(y,t)= \int_0^\infty {\cal L}(y,s)\psi(t-s) \ ds\\
  && = \int_0^\infty {\cal L}(y,s)\left(\bmt z_2(L,t-s)\\z_4(L,t-s)\emt+ {\cal D}\bmt  w_1(t-s)\\ w_2(t-s)\emt\right)\ ds
  \eean
  for some $4\times 2$ matrix valued function ${\cal L}(y,s)$.  Notice that ${\cal L}(y,s)$ depends on the location $y$ because we are trying to estimate 
  $z(y,t) $.

 Given such a  ${\cal L}(y,s)$ for each $y\in [0,L]$ we define a $4\times 4$ matrix valued function ${\cal H}(y,y_1,s)$ by the differential equation 
 \bea \label{adj_dyn}
 \frac{\partial {\cal H}}{\partial s}(y,y_1,s)&=&{\cal H}(y,y_1,s){\cal A}_1 -{\cal L}(y,s) {\cal C}(y_1) \eea
The subscript on ${\cal A}_1$ indicates that it is involves partial differentiation with respect to $y_1$ not $y$. 

  This PDE is subject to the homogeneous boundary conditions 
  \bean
  {\cal H}_{:,2}(y,0,t)=0,&&
   \frac{\partial^2 {\cal H}_{:,2}}{\partial y_1^2}(y,0,t)=0\\ 
   \frac{\partial^2 {\cal H}_{:,2}}{\partial y_1^2}(y,L,t)=0,&&
   \frac{\partial^3 {\cal H}_{:,2}}{\partial y_1^3}(y,L,t)=0\\
    \frac{\partial {\cal H}_{:,4}}{\partial y_1}(y,0,t)=0,&&  \frac{\partial {\cal H}_{:,4}}{\partial y_1}(y,L,t)=0
   \eean
     and the initial condition
  \bean
  {\cal H}(y,y_1,0)=\delta(y-y_1) I^{4\times 4}
  \eean
  
  Then
  \bean
 && \hat{z}(y,t)=\int_0^\infty \int_0^L \left(-\frac{\partial {\cal H}}{\partial s}(y,y_1,s)+\left( {\cal H}(y,y_1,s){\cal A}_1\right)\right)\\
  && \times\ z(y_1,t-s) \ dy_1\ ds
\eean 
Note the parenthesis around $\left( {\cal H}(y,y_1,s){\cal A}_1\right)$.  The parenthesis  implies that we apply the matrix partial differential operator ${\cal A}_1$ to $ {\cal H}(y,y_1,s)$ before doing anything else with this term.

Now we integrate by parts with respect to $s$ 
\bean
&&-\int_0^\infty \int_0^L \frac{\partial {\cal H}}{\partial s}(y,y_1,s)z(y_1,t-s)\ dy_1\ ds\\
&&=-\int_0^L \left[{\cal H}(y,y_1,s)z(y_1,t-s)\right]_0^\infty \ dy_1\\
& &+\int_0^\infty \int_0^L {\cal H}(y,y_1,s) \frac{\partial z}{\partial s}(y_1,t-s)\ dy_1\ ds
  \eean
We expect that the state $z(y_1,t-s)$ in the far past $s>>0$ to be  irrelevant to the estimate  $\hat{z}(y_1,t)$ of the  state at time $t$ so we assume that  ${\cal H}(y,y_1,s)z(y_1,t-s)\to 0 $ as $s\to \infty$.  Hence
 \bean
&&-\int_0^\infty  \int_0^L\frac{\partial {\cal H}}{\partial s}(y,y_1,s)z(y_1,s)\ dy_1 \ ds\\
&&= z(y,t)+\int_0^\infty \int_0^L{\cal H}(y,y_1,s)\\
&& \times\ \left(-{\cal A}_1z(y_1,t-s)\right)+{\cal B}(y_1)v(t-s)\ dy_1\ ds
\eean 
 
 Then we integrate by parts with respect to $y_1$ using the boundary conditions 
 \bean
&&-\int_0^L\left({\cal H}(y,y_1,s) {\cal A}_1\right) z(y_1,t-s)\ dy_1\\
&&=-\int_0^L{\cal H}(y,y_1,s)\left({\cal A}_1z(y_1,t-s)\right)\ dy_1
 \eean
 
 Putting this  all together yields
 \bean
  && \hat{z}(y,t)= z(y,t)+\int_0^\infty \int_0^L {\cal H}(y,y_1,s){\cal B}(y_1)v(t-s)\ dy_1\ ds\\
  && +\int_0^\infty {\cal L}(y,s) {\cal D}\bmt w_1(t-s)\\w_2(t-s)\emt
 \eean
 So the estimation error $\tilde{z}(y,t)=z(y,t)-\hat{z}(y,t)$ is given by
  \bean
   \tilde{z}(y,t)&=& -\int_0^\infty \int_0^L {\cal H}(y,y_1,s){\cal B}(y_1)v(t-s)\ dy_1\ ds\\
  && -\int_0^\infty {\cal L}(y,s) {\cal D}\bmt w_1(t-s)\\w_2(t-s)\emt \ ds
 \eean
Because $v(t)$ and $w(t)$ are standard white Gaussian noises,  the  covariance of the error is
\bea \label{critK}
&&\int_0^\infty \iint_{\cal S}{\cal H}(y,y_1,s){\cal B}(y_1){\cal B}\tr (y_2){\cal H}\tr (y,y_2,s)\ dA \ ds\nonumber \\
&&+\int_0^\infty {\cal L}(y,s) {\cal D}{\cal D}\tr  {\cal L}\tr (y,s)\ ds \label{adj_crit}
\eea

 For each $y\in [0,L]$ and for each pair of corresponding rows of ${\cal H}(y,y_1,s)$ and ${\cal L}(y,s)$ we have an LQR in adjoint form with state ${\cal H}_{i,:}(y,y_1,s)$, control ${\cal L}_{i,:}(y,s)$, linear dynamics (\ref{adj_dyn})
 and a quadratic criterion (\ref{adj_crit}).  We can leave this adjoint LQR in matrix form as  the optmal feedback gain ${\cal K}(y,y_1)$ is the same for all rows,
 \bean
  {\cal L}(y,s) &=& \int_0^L  {\cal H}(y,y_1,s){\cal K}(y,y_1)\ dy_1
 \eean
 That is not to say that the rows of $ {\cal L}(y,s) $ are all the same because the rows of $ {\cal H}(y,y_1,s) $ can be different,
 We are trying to estimate  $z(y,t)$ which is why ${\cal H}(y,y_1,s)$ and ${\cal K}(y,y_1)$ might depand on $y$.
 
For each $y\in [0,L]$ let ${\cal P}(y,y_1,y_2)$ be a $6\times 6$ continuous matrix value function of $y_1,y_2$ satisfying the homogeneous  boundary conditions
 \bean
 {\cal P}(y,y_1,0)=0,&& {\cal P}(y,y_1,L)=0\\
  {\cal P}(y,0,y_2)=0,&& {\cal P}(y,L,y_2)=0
 \eean
Using the Fundamental Theorem of Calculus assuming ${\cal H}(y,y_1,s)\to 0$ as $s\to \infty$ we have
 \bean
 0&=& \iint_{\cal S} {\cal H}(y,y_1,0){\cal P}(y,y_1,y_2){\cal H}\tr(y,y_2,0)\ dA\\
 && + \int_0^\infty {d\over ds}\iint_{\cal S}{\cal H}(y,y_1,s){\cal P}(y,y_1,y_2){\cal H}\tr(y,y_2,s)\ dA\\
 &=& \iint_{\cal S} {\cal H}(y,y_1,0){\cal P}(y,y_1,y_2){\cal H}\tr(y,y_2,0)\ dA\\
 &&  \int_0^\infty \iint_{\cal S}\left( {\cal H}(y,y_1,s)  { \cal A}_1- {\cal L}(y,s){\cal C}(y_1)  \right)  \\
 && \times   {\cal P}(y,y_1,y_2){\cal H}\tr(y,y_2,s)\ dA\\
 &&+ \int_0^\infty \iint_{\cal S}{\cal H}(y,y_1,s) {\cal P}(y,y_1,y_2)\\
 &&\times
 \left( {\cal H}(y,y_2,s)  { \cal A}_2- {\cal L}(y,s){\cal C}(y_2)  \right)\tr   \ dA
 \eean
 
 We formally integrate  by parts as if ${\cal P}(y,y_1,y_2)$ were $C^2$ in $y_1,y_2$
 \tiny
 \bea \label{zero}
 &&0= \iint_{\cal S} {\cal H}(y,y_1,0){\cal P}(y,y_1,y_2){\cal H}\tr(y,y_2,0)\ dA\\ \nonumber
 & &+ \int_0^\infty \iint_{\cal S}{\cal H}(y,y_1,s)\left({ \cal A}_1 {\cal P}(y,y_1,y_2)\right) {\cal H}\tr(y,y_2,s)\ dA\ dt\\  \nonumber
 &&+ \int_0^\infty \iint_{\cal S}{\cal H}(y,y_1,s)\left( {\cal P}(y,y_1,y_2) {\cal A}\tr_2\right) {\cal H}\tr(y,y_2,s)\ dA\ dt\\ \nonumber
 &&- \int_0^\infty \iint_{\cal S} {\cal L}(y,s){\cal C}(y_1) {\cal P}(y,y_1,y_2)  {\cal H}\tr(y,y_2,s)\ dA\ dt\\ \nonumber
  &&- \int_0^\infty \iint_{\cal S} {\cal H}(y,y_2,s){\cal P}(y,y_1,y_2) {\cal C}\tr(y_2){\cal L}\tr(y,s)  \ dA\ dt  \nonumber
 \eea
 \normalsize
 which reduces to
  \tiny
 \bea \label{zero1}
 &&0= \iint_{\cal S} {\cal H}(y,y_1,0){\cal P}(y,y_1,y_2){\cal H}\tr(y,y_2,0)\ dA\\ \nonumber
 & &+ \int_0^\infty \iint_{\cal S}{\cal H}(y,y_1,s)\left({ \cal A}_1 {\cal P}(y,y_1,y_2)\right) {\cal H}\tr(y,y_2,s)\ dA\ dt\\  \nonumber
 &&+ \int_0^\infty \iint_{\cal S}{\cal H}(y,y_1,s)\left( {\cal P}(y,y_1,y_2) {\cal A}_2\tr\right) {\cal H}\tr(y,y_2,s)\ dA\ dt\\ \nonumber
 && -\int_0^\infty\int_0^L   {\cal L}(y,s) \bmt  {\cal P}_{2,:}(y,L,y_2) \\{\cal P}_{4,:}(y,L,y_2) \emt {\cal H}\tr(y,y_2,s)          \ dy_2\ dt\\ \nonumber
 && -\int_0^\infty\int_0^L   {\cal H}(y,y_1,s)\bmt {\cal P}_{:,2}(y,y_1,L)&{\cal P}_{:,4}(y,y_1,L)\emt \nonumber \\
 && \times{\cal L}\tr(y,s)  \ dy_1\ dt  \nonumber
 \eea
  \normalsize

 We add the right side of (\ref{zero1}) to the criterion (\ref{critK}) to be minimized to get an equivalent criterion
 \tiny
 \bea  \label{eqcritK}
 &&\iint_{\cal S} {\cal H}(y,y_1,0){\cal P}(y,y_1,y_2){\cal H}\tr(y,y_2,0)\ dA\\ \nonumber
 & &+ \int_0^\infty \iint_{\cal S}{\cal H}(y,y_1,s)\left({ \cal A}_1 {\cal P}(y,y_1,y_2)\right) {\cal H}\tr(y,y_2,s)\ dA\ dt\\  \nonumber
 &&+ \int_0^\infty \iint_{\cal S}{\cal H}(y,y_1,s)\left( {\cal P}(y,y_1,y_2) {\cal A}\tr_2\right) {\cal H}\tr(y,y_2,s)\ dA\ dt\\ \nonumber
 & &+ \int_0^\infty \iint_{\cal S}{\cal H}(y,y_1,s){\cal B}(y_1){\cal B}\tr(y_2){\cal H}\tr(y,y_2,s)\ dA\ dt\\  \nonumber
  && -\int_0^\infty\int_0^L   {\cal L}(y,s) \bmt  {\cal P}_{2,:}(y,L,y_2) \\{\cal P}_{4,:}(y,L,y_2) \emt {\cal H}\tr(y,y_2,s)          \ dy_2\ dt\\ \nonumber
 && -\int_0^\infty\int_0^L   {\cal H}(y,y_1,s)\bmt {\cal P}_{:,2}(y,y_1,L)&{\cal P}_{:,4}(y,y_1,L)\emt \nonumber \\
 && \times{\cal L}\tr(y,s)  \ dy_1\ dt  \nonumber+ \int_0^\infty {\cal L}(y,s){\cal D}{\cal D}\tr{\cal L}\tr(y,s)
  \eea
  \normalsize

 We seek a ${\cal K}(y,y_1)$ such that the time integrand in (\ref{eqcritK}) is a perfect square of the form
 \bea \label{psK}
 &&\left({\cal L}(y,s)-\int_0^L {\cal H}(y,y_1,s){\cal K}(y,y_1)\ d y_1\right){\cal D}{\cal D}\tr\\
  && \times \nonumber
  \left({\cal L}(y,s)-\int_0^L {\cal H}(y,y_2,s){\cal K}(y,y_2)\ d y_2\right)\tr
 \eea
 The terms quadratic in ${\cal L}(y,s)$ match so we equate terms bilinear
  in ${\cal L}(y,s)$ and ${\cal H}(y,y_2,s)$ and obtain
 \bean
 {\cal D}{\cal D}\tr{\cal K}\tr (y,y_2)&=& \bmt {\cal P}_{2,:}(y,L,y_2)\\  {\cal P}_{4,:}(y,L,y_2)\emt
 \eean
 so we set 
 \bea \label{calK}
 {\cal K}(y,y_2)= \bmt {\cal P}_{:,2}(y,y_1,L)&  {\cal P}_{:,4}(y,y_1,L)\emt \left( {\cal D}{\cal D}\tr\right)^{-1}
 \eea
 
 Finally we match terms bilinear in  ${\cal H}(y,y_1,s)$ and ${\cal H}(y,y_2,s)$ which yields the Riccati PDE of the 
 infinite dimensional Kalman filter
 \bea
 &&\left({ \cal A}_1 {\cal P}(y,y_1,y_2)\right) +\left( {\cal P}(y,y_1,y_2) {\cal A}\tr_2\right) +{\cal B}{\cal B}\tr\\
 &&= \bmt {\cal P}_{:,2}(y,y_1,L)&  {\cal P}_{:,4}(y,y_1,L)\emt \left( {\cal D}{\cal D}\tr\right)^{-1}\nonumber \\ \nonumber 
 && \times \bmt {\cal P}_{2,:}(y,L,y_2)\\  {\cal P}_{4,:}(y,L,y_2)\emt
 \eea
 Actually this is a family of Riccati PDEs paramterized by $y$,

 The optimal estimate is given by
\bean
  \hat{z}(y,t)&=& \int_0^\infty {\cal L}(y,s)\psi(t-s) \ ds \\
  &=& \int_0^\infty  \int_0^L  {\cal H}(y,y_1,s){\cal K}(y,y_1)\ dy_1 \psi(t-s) \ ds\\
    \eean 
    We make the substitution $s=t-\tau$ then 
    \bean
  \hat{z}(y,t)&=&\int_{-\infty}^t  \int_0^L  {\cal H}(y,y_1,t-\tau){\cal K}(y,y_1)\ dy_1\ \psi(\tau) \ d\tau
    \eean

 We differentiate this with respect $t$ and obtain 
 \bean
&& \frac{\partial   \hat{z}}{\partial t}(y,t)  = \int_0^L  {\cal H}(y,y_1,0){\cal K}(y,y_1)\ dy_1\ \psi(t)\\
 &&+ \int_{-\infty}^t  \int_0^L \frac{\partial  {\cal H}}{\partial t}(y,y_1,t-\tau){\cal K}(y,y_1)\ dy_1 \psi(\tau) \ d\tau\\
 \\
 &&\frac{\partial   \hat{z}}{\partial t}(y,t) ={\cal K}(y,y)\psi(t)\\
&&+\int_{-\infty}^t \int_0^L \left({\cal H}(y,y_1,t-\tau)D_1-  {\cal L}(y_1,s){\cal C}(y_1)\right)\\
&&\times {\cal K}(y,y_1)\ dy_1 \psi(\tau) \ d\tau
 \eean
 
 Then
  \bean \label{adj_dyncl}
&& \frac{\partial {\cal H}}{\partial s}(y,y_1,s)={\cal H}(y,y_1,s){\cal A}\\
&& \nonumber -\int_0^L{\cal H}(y,y_1,s){\cal K}(y,y_1) {\cal C}(y_1)\ dy_1\\
 && ={\cal H}(y,y_1,s)  {\cal F}(y,y_1)
 \eean
where ${\cal F}(y,y_1)$ is the differential operator 
 \bean
 {\cal F}(y,y_1)z(y_1,t)&=&{\cal A}_1z(y_1,t)- {\cal K}(y,y_1){\cal C}(y_1)
  \eean

 Put another way 
 \bean \label{adj_dyncl}
 \frac{\partial {\cal H}}{\partial s}(y,y_1,s)&=& {\cal H}(y,y_1,s){\cal F}(y,y_1)
 \eean
 and so for each $y,y_1$
 \bean
   {\cal H}(y,y_1,\tau) {\cal H}(y,y_1,t-\tau)&=&{\cal H}(y,y_1,t)
 \eean
 We differentiate this with respect  to $\tau$ and obtain
  \bean
0&=&{\cal H}(y,y_1,\tau){\cal F}(y,y_1){\cal H}(y,y_1,t-\tau)\\
&&+{\cal H}(y,y_1,\tau)  \frac{\partial   {\cal H}}{\partial \tau }(y,y_1,t-\tau)\\
 \eean
 so we conclude that  
 \bean
 \frac{\partial   {\cal H}}{\partial \tau }(y,y_1,t-\tau)&=&{\cal F}(y,y_1){\cal H}(y,y_1,t-\tau)
 \eean
 and
 \bean
 \frac{\partial   {\cal H}}{\partial t }(y,y_1,t-\tau)&=&-{\cal F}(y,y_1){\cal H}(y,y_1,t-\tau)
 \eean

  We differentiate the optimal estimate with respect $t$ and obtain 
 \bean
&&  \frac{\partial   \hat{z}}{\partial t}(y,t)   ={\cal K}(y,y)\psi(t)\\
&&+ \int_{-\infty}\tr  \int_0^L \frac{\partial  {\cal H}}{\partial t}(y,y_1,t-\tau){\cal K}(y,y_1)\ dy_1 \psi(\tau) \ d\tau\\
&&={\cal K}(y,y)\psi(t)\\
&&+ \int_{-\infty}\tr  {\cal F}(y)\int_0^L {\cal H}(y,y_1,t-\tau){\cal K}(y,y_1)\ dy_1 \psi(\tau) \ d\tau\\
\eean
so we conclude that
\bean
  \frac{\partial   \hat{z}}{\partial t}(y,t)  = {\cal F}(y) \hat{z}(y,t) +{\cal K}(y,y)\psi(t)
\eean
or
\bea 
  \frac{\partial   \hat{z}}{\partial t}(y,t)  = {\cal A}(y) \hat{z}(y,t) +{\cal K}(y,y)\left( {\psi}(t)-\bmt 
 \hat{z}_2(L,t)\\ \hat{z}_4(L,t)\emt \right)
 \eea
 The boundary conditions on  $ \hat{z}(y,t)$ are the same as those (\ref{zzbc}) on ${z}(y,t)$,
 \bea \label{zhatbc}
 \hat{z}_1(0,t)=0, && \frac{\partial^2 \hat{z}_1}{\partial  y^2}(0,t)=B_1u_(t)\\ \nonumber
 \frac{\partial^2 \hat{z}_1}{\partial  y^2}(L,t)=0,&& \frac{\partial^3 \hat{z}_1}{\partial  y^3}(L,t)=0\\ \nonumber
  \frac{\partial \hat{z}_3}{\partial  y}(0,t)=B_2u_2(t),&& \frac{\partial \hat{z}_3}{\partial  y}(L,t)=0
 \eea
 
The so called innovations process is
\bean
\tilde{\psi}(t)&=& \psi(t)-\hat{\psi}(t)
\eean
where
 \bean
 \hat{\psi}(t)&=&\bmt \hat{z}_2(L,t)\\ \hat{z}_4(L,t)\emt
 \eean
 This is the difference between the actual observation $\psi(t)$ and what we thing it should
 be $\hat{\psi}(t)$ given our optimal estimate of the state.

We can express the optimal filter as a copy of the original dynamics driven by the innovations
 \bea \label{zhatdyn}
 \frac{\partial   \hat{z}}{\partial t}(y,t)    &=& {\cal A}(y) \hat{z}(y,t) +{\cal K}(y,y) \tilde{\psi}(t)
  \eea
  subject to the boundary conditions (\ref{zhatbc}).

  If we use the estimate $\hat{z}(y,t)$ in place of the full state $z(y,t)$ in the feedback law we found before then the control input is
  \bean
  u(t)&=& \int_0^L K(y) \hat{z}(y,t)\ dy
  \eean
   We insert $u(t) $ into the plant boundary conditions (\ref{zbc}).  We also insert   it into the boundary  conditions  (\ref{zhatbc}) of the Kalman filter.  This is dynamic compensation.

   The error  $ \tilde{z}(y,t)=z(y,t)-\hat{z}(y,t)$ dynamics is given by
  \bean
  \frac{\partial   \tilde{z}}{\partial t}(y,t)  &=& {\cal F}(y) \tilde{z}(y,t) 
  \eean
subject to the homogeneous boundary conditions 
 \bea \label{ztiltbc}
 \tilde{z}_1(0,t)=0, && \frac{\partial^2\tilde{z}_1}{\partial  y^2}(0,t)=0\\ \nonumber
 \frac{\partial^2 \tilde{z}_1}{\partial  y^2}(L,t)=0,&& \frac{\partial^3 \tilde{z}_1}{\partial  y^3}(L,t)=0\\ \nonumber
  \frac{\partial \tilde{z}_3}{\partial  y}(0,t)=0,&& \frac{\partial \tilde{z}_3}{\partial  y}(L,t)=0
 \eea

  Notice that the error dynamics and its boundary conditions  only depend  on the error $\tilde{z}(y,t)$ so   we can express the combined plant and Kalman filter in the coordinates 
  $z(y,t)$ and $\tilde{z}(y,t)$.  In these coordinates the combined system is upper triangular so
if the full state feedback aymptotically stabilizes the system and if the error dynamics of the Kalman filter is asymptotically stable
then the dynamic compensator stabilizes the system.

  If the coefficient of the driving noise does not depend on $y$, ${\cal B}(y)={\cal B}$ then the family of LQRs 
  does not depend on $y$.  Therefore ${\cal L}(y,t)={\cal L}(t)$,  ${\cal H}(y,y_1,t)={\cal H}(y_1,t)$,
 ${\cal K}(y,y_1)={\cal K}(y_1)$ and ${\cal P}(y,y_1,y_2)={\cal P}(y_1,y_2)$.
  The optimal filter is given by
  \bea 
  \frac{\partial   \hat{z}}{\partial t}(y,t) (y,t)  &=& {\cal A}(y) \hat{z}(y,t) +{\cal K}(y) \tilde{\psi}(t)\\
  \eea

\section{Filter Example}
\begin{figure}
 \centering
\includegraphics[width=3.5in]{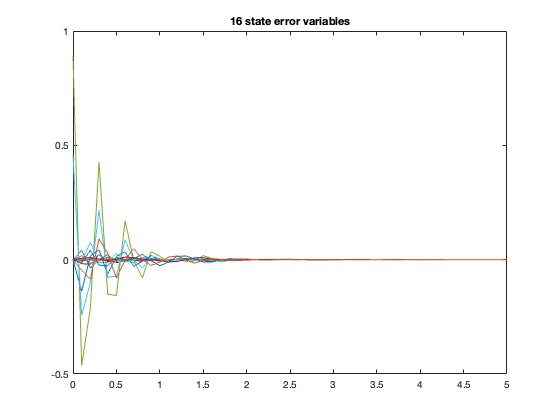}
\caption{Sixteen Dimensional Errors}  
\label{16d}      
 \end{figure}

We apply the Kalman filter to the $16^{th}$ dimensional example that we treated above.   The $16$ error variables
are shown in Figure \ref{16d}.
    
 \begin{figure}
 \centering
\includegraphics[width=3.5in]{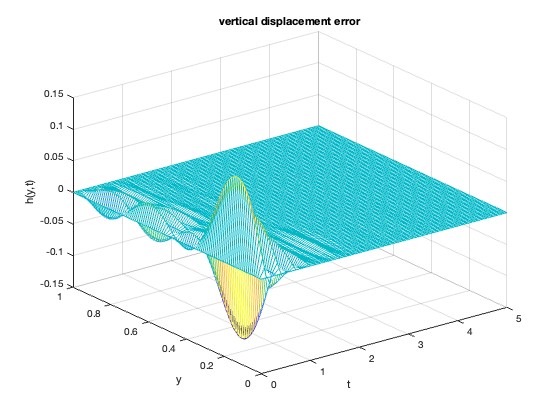}
\caption{Vertical Displacement Estimation Error}    
\label{vde}  
 \end{figure}
\begin{figure}
 \centering
\includegraphics[width=3.5in]{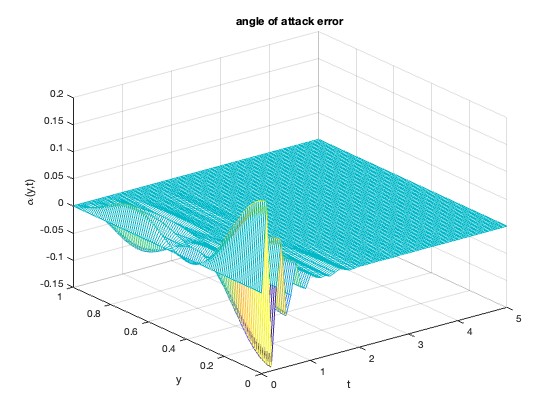}
\caption{Angle of Attack Estimation Error}       
\label{aoaerr}
 \end{figure}

 The vertical displacement error is shown in Figure \ref{vde} and the rotational displacement error is shown in Figure \ref{aoaerr} 
 Notice  that the estimation errors are smaller near the sensors at the tip, $ y=L$.

\section{ Aerodynamic Model}
The next step is to add a model of the  aerodynamical forces that are generated by the bending and twisting beam to get a model of a wing.The classical aerodynamical models are those of Wagner and Theodoresen.   Wagner's model is in the time
domain and Theodoresen's is in the frequency  domain.
They are related by the Laplace transform.

These  models are not exactly realizable by a linear time invariant system.    But over the years 	considerable effort has gone into
finding linear time invariant systems whose impulse response or transfer function approximates Wagner's or Theodoresn's model.
We shall  use the approximation to Wagner's model developed by R. T. Jones in 1938 as reported in 
\cite{BR13}.
Wagner's model is valid for an airfoil, a wing section.
We shall extend it to  model  a wing by introducing spanwise dependence.

The state space realization of Wagner's model  adds two additional states to
our four dimensional model.    
Of course each of these six state variables is a function of $ y$
and  $ t$.
  
  Aerodynamical models are parameterized by  the free stream air velocity $ U_\infty$ which we could
treat as a static state.  But then the model would become  nonlinear so we do not do so..
     
   Here is R.T. Jones' approximate state space  realization  of  Wagner's impulse  response.
   The aerodynamics state is $\xi=[\xi_1(y,t),\xi_2(y,t)]\tr $ and its dynamics is
\bean
\frac{\partial \xi}{\partial t}(y,t)&=& A\xi(y,t)+B\nu(y,t)
\eean
where
\bean
A=\bmt 0&1\\-0.0137&-0.3455\emt,&& B=\bmt 0\\1\emt
\eean
 The aerodynamic input is
\bean
\nu(y,t)&=&-\dot{h}(y,t)+\left(U_\infty + b\left({1\over 2}-a\right)\right)\alpha(y,t)
\eean
The aerodynamics output is 
\bean
\psi(y,t)&=& C\xi(y,t)+D\nu(y,t)
\eean
where
\bea
C=\bmt 0.0068&0.1080\emt,&& D=\bmt0.5\emt
\eea

The lift per unit span is 
\bean
L&=&\rho_\infty\pi b\left(-\ddot{h}+U_\infty\dot\alpha-ba \ddot{\alpha}\right)\\
&&+\pi  \rho_\infty U_\infty\left(-\dot{h}+U_\infty\alpha+b\left({1\over 2}-a\right)\dot{\alpha}\right)\\
&& -2\pi b  \rho_\infty U_\infty\psi\
\eean

The moment per unit span is 
\bean
M_{ea}&=&\rho_\infty\pi b^2\left(-a\ddot{h}+\left(a-{1\over 2}\right) U_\infty\dot{\alpha}-b\left(a^2+{1\over 8}\right)\ddot{\alpha}\right)\\
 &&\pi b \rho_\infty U_\infty\left(a+{1 \over 2}\right)\left(-\dot{h}+U_\infty\alpha+b\left({1\over2}-a\right)\dot{\alpha}\right)\\
 &&-2\pi b \rho_\infty U_\infty\left(a+{1\over 2}\right)\psi 
\eean
where $ \psi(y,t) $   is the aerodynamic output.

Typical parameter values are 
\bean
\bmt
L&\mbox{span} & 15 m\\
b&\mbox{half chord} &0.5m\\
a&\mbox{elastic axis relative to mid chord} &0.0m\\
\rho_\infty& \mbox{freestream air density}&0.0889 kg/m^3\\
U_\infty&  \mbox{freestream velocity}& 45m/s\\
\emt
\eean
The model is  from \cite{ HMA19} and the 
parameters are from \cite{BR13}.

We conjectured that the combination of  $\Phi_n(y), \ n=1,2,3,\ldots$ and $\Theta_m(y),\ m=0.1,2,\ldots$ constitutes a Riesz basis for $L^2[0,L]$.
     We expand the aerodynamic state in these functions
\bean
\bmt \sigma_1(y,t)\\ \sigma_2(y,t)\emt&=&\sum_{n=1}^\infty\bmt \xi_{1,n}(t)\\ 
\xi_{2,n}(t)\emt \Phi_n(y)\\&&+\sum_{m=0}^\infty\bmt \xi_{3,m}(t)\\ \xi_{4,m}(t)\emt\Theta_{m}(y)
\eean
We differentiate with respect to time
\bean
\bmt \dot{\sigma_1}(y,t)\\ \dot{\sigma_2}(y,t)\emt&=&\sum_{n=1}^\infty\bmt \dot{\xi}_{1,n}(t)\\ \dot{\xi}_{2,n}(t)\emt \Phi_n(y)+\sum_{m=0}^\infty\bmt \dot{\xi}_{3,m}(t)\\ \dot{\xi}_{4,m}(t)\emt\Theta_{m}(y)
\eean

On the other hand
\small
\bean
A\bmt \sigma_1(y,t)\\ \sigma_2(y,t)\emt&=&A\sum_{n=1}^\infty\bmt \xi_{1,n}(t)\\ \xi_{2,n}(t)\emt \Phi_n(y)+A\sum_{n=0}^\infty\bmt \xi_{3,n}(t)\\ \xi_{4,n}(t)\emt\Theta_{n}(y) 
\eean
\bean
Bw(t)&=&B\left(-  \dot{h}(y,t)+\left(U_\infty + b\left({1\over 2}-a\right)\right)\alpha(y,t)\right)\\
&=&-B  \sum_{n=1}^\infty \zeta_{2,n}(t) \Phi_n(y) \\
&&+B\left(U_\infty+ b\left({1\over 2}-a\right)\right)\sum_{n=0}^{\infty}\zeta_{3,n}(t)\Theta_n(y)
\eean 

We multiply by the elements of the dual basis and get the modal dynamics
\bean
\bmt \dot{\xi}_{1,n}(t)\\ \dot{\xi}_{2,n}(t)\emt &=& A\bmt \xi_{1,n}(t)\\ \xi_{2,n}(t)\emt-B\zeta_{2,n}\\
\bmt \dot{\xi}_{3,n}(t)\\ \dot{\xi}_{4,n}(t)\emt &=& A\bmt \xi_{3,n}(t)\\ \xi_{4,n}(t)\emt+B\left(U_\infty+ b\left({1\over 2}-a\right)\right)\zeta_{3,n}(t)
\eean

The modal decomposition of the aerodynamic output is 
\small
\bean
&&y_a(y,t)= C\sum_{n=1}^\infty\bmt \xi_{1,n}(t)\\ \xi_{2,n}(t)\emt \Phi_n(y)+C\sum_{n=0}^\infty\bmt \xi_{3,n}(t)\\ \xi_{4,n}(t)\emt\Theta_{n}(y) \\
&&-D\sum_{n=1}^\infty \zeta_{2,n}\Phi_n(y)+D\left(U_\infty+ b\left({1\over 2}-a\right)\right)\sum_{n=0}^\infty\zeta_{3,n}(t)\Theta_{n}(y)
\eean

Our approximate model of beam has a sixteen dimensional state
\bean
z_1(y,t)&\approx&\bmt \Phi_1(y)&\ldots & \Phi_4(y)\emt \zeta_1(t) \\
z_2(y,t)&\approx&\bmt \Phi_1(y)&\ldots & \Phi_4(y)\emt \zeta_2(t) \\
z_3(y,t)&\approx&\bmt \Theta_0(y)&\ldots & \Theta_3(y)\emt \zeta_3(t) \\
z_4(y,t)&\approx&\bmt \Theta_0(y)&\ldots & \Theta_3(y)\emt \zeta_4(t) \
\eean 
where $ \zeta_i(t)$ is  $4\times 1$.
The approximate model of aerodynamics is also sixteen dimensional.
\bean
\sigma_1(y,t)&\approx & \bmt \Phi_1(y)&\ldots & \Phi_4(y)\emt  \xi_1(t)\\
&&+\bmt \Theta_0(y)&\ldots & \Theta_3(y)\emt  \xi_3(t)\\
\sigma_2(y,t)&\approx & \bmt \Phi_1(y)&\ldots & \Phi_4(y)\emt  \xi_2(t)\\
&&+\bmt \Theta_0(y)&\ldots & \Theta_3(y)\emt  \xi_4(t)
\eean
where $ \xi_i(t)$ is $ 4\times 1$.
So the combined approximate model is $ 32$ dimensional.

\begin{figure}
 \centering
\includegraphics[width=4in]{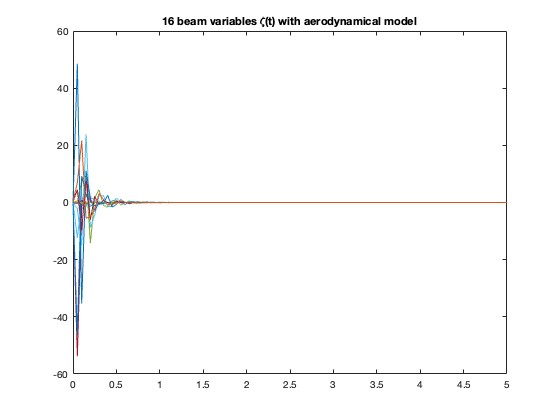}
\caption{Beam Variables  $\zeta(t)$ with Aero Model}  
\label{BVAM}      
 \end{figure}

Figure \ref{BVAM} shows the $16$ beam variables of $32$ variable combined model being stabilized by the $32$ dimensional full state feedback.

\begin{figure}
 \centering
\includegraphics[width=4in]{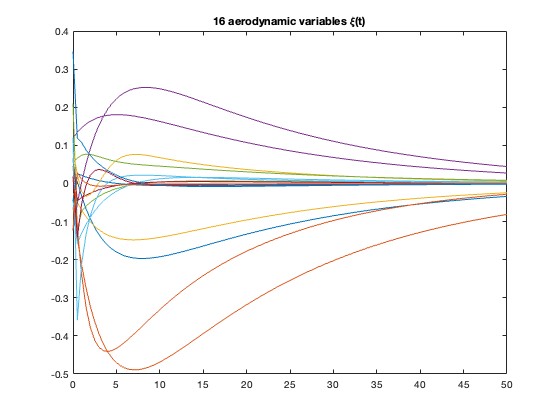}
\caption{Aero Variables  $\xi(t)$ }      
\label{AVAM}
 \end{figure}

Figure \ref{AVAM} shows the $16$ aerodynamic  variables of $32$ variable combined model being stabilized by the $32$ dimensional full state feedback.
Notice how slow the stabilization is.   The bending and twisting controls have very limited authority if any over the aerodynamical states.  We can still use LQR because the eigenvalues  of $ A$ are
$
-0.0457$ and $
   -0.2998
$
so aero model is lightly damped.

\begin{figure}
 \centering
\includegraphics[width=4in]{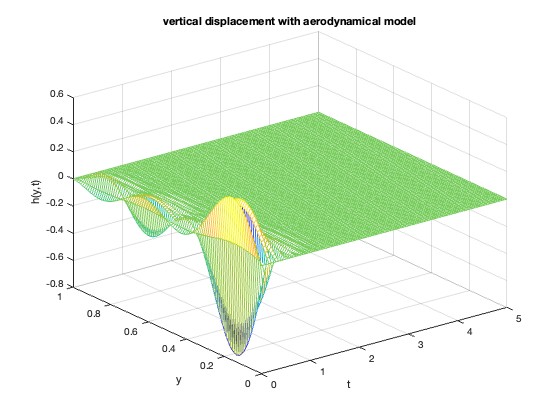}
\caption{Vertical Displacement with Aero Model}    
\label{VDAM}
 \end{figure}

Figure \ref{VDAM} shows the  $32$ dimensional full state feedback stabilizing the vertical displacement.

\begin{figure}
 \centering
\includegraphics[width=4in]{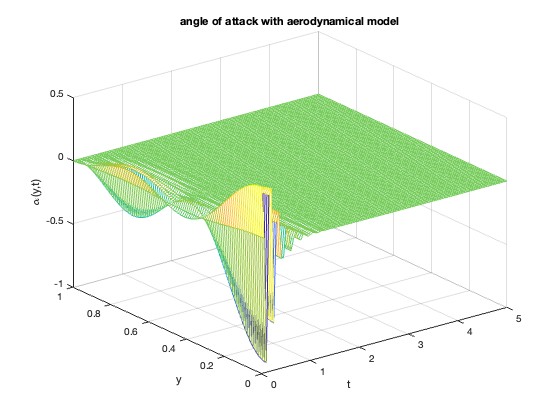}
\caption{Angle of Attack  with Aero Model}    
\label{AoAAM} 
 \end{figure}

 Figure \ref{AoAAM} shows the  $32$ dimensional full state feedback stabilizing the rotational  displacement.

 The combined beam and aero model is $ 32$ dimensional.
 The closed loop spectrum consists of $16$ complex eigenvalues
 and  $ 16 $ real eigenvalues.  
 The complex eigenvalues are loosely associated to the closed loop dynamics of the beam
 while the real eigenvalues are loosely associated to the closed loop dynamics  of the aero model.
  The real eigenvalues are very lightly damped, the largest one is $ -0.04569$. 

\section{Conclusion} 
We started by considering the bending and  twisting of a beam and showed how its oscillations can be damped by an LQR 
derived full state feedback using two point actuators located at the base of the beam.   Then we derived a Kalman filter
that processed two point measurements at the tip of the beam to obtain an estimate of the full state of the beam. This yields an Linear
Quadratic Gaussian (LQG) synthesis.
By adding a model of the aerodynamic forces generated by the bending and  twisting of the beam we converted the model
of the beam into a model of a wing and we showed that the LQR full state feedback stabilized the model of the wing.


\newpage

\begin{thebibliography}{DDDD}

\bibitem{BS97}
H.~T.~ Banks, R.~C.~Smith, D.~E.~Brown, R.~J.~Silcox and V.~L.~Metcalf,
Experimental Confirmation of a PDE-Based Approach to the Design of Feedback Controls.
SIAM J. Control and  Optimization, 
Vol. 35, No. 4, pp. 1263--1296, July 1997

\bibitem{BAH96}
R.~L,~Bisplinghoff, H.~Ashley and R.~Halfman,
{\it Aeroelasticity}, Dover, 1996.

\bibitem{BR13}
S.~L.~Brunton and C.~W.~Rowley,
Empirical state-space representations for Theodorsen's lift model, J. of Fluids and Structures, V. 38, pp.174-186, 2013


\bibitem{Do94}
E.~H.~Dowell, ed.,
{\it A Modern Course in Aeroelasticity}, Kluwer Academic Publishers, Dordrecht, 1994.

\bibitem{EBB78}
J.~W.~Edwards, J.~V.~ Breakwell, A.~E.~Bryson.,  Active flutter control using generalized unsteady aerodynamic theory. Journal of
Guidance and Control 1, 32-40, 1978.

\bibitem{HMA19}
A.~Hossein Modaress-Aval, F.~Bakhtiari-Nejad, E.~H.~Dowell, D.~Peters, and H.~Shahverdi 
A comparative study of nonlinear aeroelastic models for high
aspect ratio wings,
J. of Fluids and Structures, v. 85, pp. 249-274, 2019


\bibitem{HMA23}
A.~Hossein Modaress-Aval, F.~Bakhtiari-Nejad, E.~H.~Dowell, H.~Shahverdi and D.~Peters,
Comparative Study of Beam and Plate Theories for Moderate Aspect Ratio Wings,
AIAA Journal, v. 61, pp. 859-874, 2023.

\bibitem{Kr25}
A.~J.~ Krener,  Boundary Stabilization of a Bending and Twisting Beam by Linear 
       Quadratic Regulation, Proceedings of the 2025 ACC, Denver.



\bibitem{KM18}
K.~Menon and R.~Mittal. Computational Modelling and Analysis of Aeroelastic Flutter, 2018 Fluid
Dynamics Conference, AIAA AVIATION Forum, (AIAA 2018-3080)

\bibitem{P08}
D.~A.~Peters,  Two-dimensional incompressible unsteady airfoil theory, an overview, Journal of Fluids and Structures 24, pp. 295-312,
2008.






\end{thebibliography}
\end{document}